\documentstyle [amssymb,twoside]{amsart}
\newtheorem{thm}{Theorem}[section]
\newtheorem{prop}[thm]{Proposition}
\newtheorem{cor}[thm]{Corollary}
\newtheorem{lem}[thm]{Lemma}
\newtheorem{remark}{Remark} 
\newtheorem{dfn}{Definition} 

\def\b0{{\bf 0}}

\def\bI{{\bf I}}

\def\b1{{\bf 1}}

\def\bR{{\bf R}}

\def\bZ{{\bf Z}}
\def\bC{{\bf C}}

\def\Blat{\mbox{\it \raise2pt\hbox{"}\kern-2pt H}}
\def\pvup{\rlap{\ ${}^{p\atop{\hbox{${}^{\vee}$}}}$}\cdots}
\def\rvup{\rlap{\ ${}^{r\atop{\hbox{${}^{\vee}$}}}$}\cdots}
\def\qvup{\rlap{\ ${}^{q\atop{\hbox{${}^{\vee}$}}}$}\cdots}
\begin{document}
\begin{center}
{\center{\Large{\bf
Logarithmic vector fields and multiplication table.
} }}
\end{center}
\begin{center}
 \vspace{2pc}
\begin{flushright}
{\it Dedicated to the 61st birthday of Kyoji Saito}
\end{flushright}
 \vspace{2pc}
{ \center{\large{ Susumu TANAB\'E }}} \noindent

 \begin{minipage}[t]{10.2cm}
{\sc Abstract.} {\em This is a review article on the
Gauss-Manin system associated to the complete intersection
singularities of projection. We show how the
logarithmic vector fields appear as coefficients to  the
Gauss-Manin system (Theorem ~\ref{thm3}).
 We examine further how the multiplication table
on the Jacobian quotient module calculates the
logarithmic vector fields tangent to
the discriminant and the bifurcation set (Proposition ~\ref{prop31}, 
Proposition~\ref{prop52}).
As applications, we establish signature formulae for Euler characteristics of
real hypersurfaces (Theorem ~\ref{thm42})
and real complete intersections (Theorem ~\ref{thm51})
by means of these fields.}
\end{minipage} \hfill
\end{center}
 \vspace{1pc}

{ \center{\section{Introduction}} }
This is a review article on the
Gauss-Manin system associated to the isolated complete intersection
singularities (i.c.i.s.)
of projection and objects tightly related with them.
The notion of i.c.i.s. of projection has been picked up among
general i.c.i.s. by Viktor Goryunov
\cite{Gory1}, \cite{Gory2} as good models to which
many  arguments on the hypersurface singularities can be applied
(see for example Theorem ~\ref{thm1}, Lemma ~\ref{lem25}).
All isolated hypersurface singularities can be
considered as a special case of the
i.c.i.s.  
of projection.
Many of important quasihomogeneous i.c.i.s.
are also i.c.i.s. of projection.

The main aim of this article is to transmit 
the message that the multiplication tables defined on
different quotient rings calculate important
data both on analytic and topological characterisation of the 
i.c.i.s. of projection.
We show that the multiplication table
on the Jacobian quotient module in $({\mathcal O}_{\tilde X \times S})^k$ 
calculates the logarithmic vector fields (i.e. the coefficients to the
Gauss-Manin system defined for the period integrals) tangent to
the discriminant and the bifurcation set (Proposition ~\ref{prop31}, 
Proposition~\ref{prop52}) of the i.c.i.s. of projection.
This idea is present already in the works by Kyoji Saito \cite{Saito3}
and James William Bruce \cite{bru} for the case of hypersurface singularities
(i.e. $k=1$).

On the other hand, as applications, 
we establish signature formulae for Euler characteristics of
real hypersurfaces (Theorem ~\ref{thm42})
and real complete intersections (Theorem ~\ref{thm51})
by means of logarithmic vector fields. These are paraphrase
of results established by Zbigniew Szafraniec \cite{Sz1}.


 \vspace{1.5pc}
\footnoterule

\footnotesize{AMS Subject Classification: 14M10 (primary), 32B10, 14P05
(secondary).

Key words and phrases: complete intersections,
Gauss-Manin system, real algebraic sets.}

 \footnotesize{Partially supported by ICTP (Trieste), JSPS grant in aid 
of Prof. Toru Ohmoto (Hokkaido Univ.).}
\normalsize

\newpage

It is well known in the study of real algebraic geometry,
Oleg Viro's patchworking method (\cite{Viro})
furnishes us with a relatively simple and effective method
to construct various nonsingular real plane projective 
algebraic curves of a given degree $m$ with different isotopy types.
As this method is based on perturbations 
of singular curves with quasihomogeneous singularities,
our study on the versal deformation of hypersurface 
singularities fits into the context of real algebraic geometry. 
We shall notice that Viro's patchworking method does not describe
all possible curves corresponding to the full deformation parameter values
outside the real discriminant. 

The deformation
parameter values $s \in \bR^\mu$
that can be treated by Viro's method are
located (on a quasihomogeneous curve) 
in certain specially selected real components of the complement
to the discriminant. This situation is explained by  the 
essential use of regular triangulation of the
Newton polyhedron of the defining equation 
$F(x,s)$ in his construction.
At the end of \S 6, Example 2, we indicate cases of real curves with
different Euler characteristics that are impossible to distinguish
after patchworking method. We hope that this approach would
give a new complementary tool to the topological 
study of real algebraic curves.

The author expresses his gratitude to Aleksandr Esterov
who drew his attention to the utility of multiplication table
and proposed the first version of Theorem~\ref{thm51}. The 
main part of this work has been accomplished during 
author's stay at the International Centre for Theoretical Physics
(Trieste) and Hokkaido University 
where the author enjoyed fruitful working condition. 
The author expresses his deep gratitude to the concerned institutions
and to Prof.Toru Ohmoto who gave him an occasion to 
report part of results at RIMS (Kyoto) conference.

{ \center{\section{Complete intersection of projection}} }

Let us consider a $k-$tuple of holomorphic germs
$$\vec f(x,u) = (f_1(x,u), \cdots, f_k(x,u)) \in ({\mathcal O}_{X})^k
\leqno(2.1)$$
in the neighbourhood of the origin for $X=(\bC^{n+1},0).$
This is a 1- parameter deformation of the germ
$$\vec f^{(0)}(x)= (f_1(x,0), \cdots, f_k(x,0)) \in ({\mathcal O}_{\tilde
X})^k \leqno(2.2)$$
for $\tilde X=(\bC^{n},0).$

After \cite{Gory1} we introduce the notion of
$R_+$equivalence of projection.
Let $p:  \bC^{n+1} \longrightarrow \bC$
be a nondegenerate linear projection i.e. $dp \not =0$.

\begin{dfn} We call the diagram
$$ Y \hookrightarrow \bC^{n+1} \longrightarrow^p \bC,$$
 the projection of the variety $Y \hookrightarrow \bC^{n+1}$
on the line.
Two varieties $Y_1, Y_2$ belong to the same
$R_+$ equivalence class of projection if there exists a biholomorphic
mapping from $\bC^{n+1}$ to $\bC^{n+1}$ that preseves the projection
and induces a translation $p \rightarrow p+const$ on the line.
\end{dfn}
In this way,  we are led to the definition of an equivalence class up to the 
following ideals,
$$ T_{f}= {\mathcal O}_{X}
\langle \frac{\partial \vec f}{\partial x_1}, \cdots,
\frac{\partial \vec f}{\partial x_n}\rangle +
\vec f^\ast\left({m}_{\bC^k,0}\right)
\cdot ({\mathcal O}_{X})^k \leqno(2.3)$$
and
$$T_f^+ := T_f + \bC \frac{\partial \vec f}{\partial u} \leqno(2.4)$$
that is nothing but the tangent space to the germ of $R_+$
equivalence class of projection.
We introduce the spaces
$$Q_f:= ({\mathcal O}_{X})^k/ T_f, \leqno(2.5)$$
$$Q_f^+:= ({\mathcal O}_{X})^k/ T_f^+. \leqno(2.6)$$
We remark that though $T_f^+ $ is not necessarily an ideal the quotien
$Q_f^+$ can make sense.
Assume that $Q_f$ is a finite dimensional $\bC$ vector space.
In this case,
we call the number $ \tau:=dim_{\bC} Q_f^+$
the  $R_+-$ codimension  of projection
 $ \mu:=dim_{\bC} Q_f$ the multiplicity of the critical point
$(x,u)=0$ of the height function $u$ on 
 $X_0:= \{ (x,u) \in X;
f_1(x,u)=\cdots= f_k(x,u)=0\}$.
We denote by $\langle \vec e_1(x,u), \cdots,  \vec e_\tau(x,u) \rangle $
the basis of the $\bC$-vector space $Q_f^+$.
If $\tau < \infty$, it is easy to see that $\vec f(x,u)=0$ (resp. $\vec f(x,0)=0$)
has isolated singularity at $0 \in X$ (resp. $0 \in \tilde X$).
Let us consider a $R_+$- versal  deformation of $\vec f^{(0)}(x)$
$$\vec F(x,u,t) =\vec f^{(0)}(x)+\vec e_0(x,u) + t_1 
\vec e_1(x,u)+ \cdots+ t_\tau \vec e_\tau(x,u), \leqno(2.7)$$
with $\vec e_0(x,u) = \vec f(x,u) - \vec f(x,0).$
We consider the deformation of $X_0$ as follows
$$X_t:= \{(x,u) \in X; \vec F(x,u,t) = \vec 0 \}, \leqno(2.8)$$
that is also a $(\tau+1)$-dimensional deformation of the germ
$\tilde X_0:= \{ x \in \tilde X;f_1(x,0)=\cdots= f_k(x,0)=0\}.$
The following fact is crucial for further arguments.
\begin{thm}(\cite{Gory1}, Theorem 2.1)
For the k-tuple of holomorphic germs $(2.1)$
with $0 < \mu <+\infty$, we have the equality $\mu = \tau+1$.
\label{thm1}
\end{thm}

Recently a conceptual understanding in terms of homological algebra
of this phenomenon appeared. See \cite{degre}, \S 3.

Futher,
in view of the Theorem ~\ref{thm1} we make use of the notation,
$S=(\bC^{\tau+1},0)=(\bC^\mu,0)$,$s$
$=$ $(u,t) \in S,$ $s_0=u, s_i=t_i,  1 \leq i \leq \tau.$
We will denote the deformation
parameter space $t \in T= (\bC^\tau,0)$.

Let $I_{C_0} \subset {\mathcal O}_X$ be
the ideal generated by $k \times k $ minors of the marix
$ (\frac{\partial \vec f(x,u)}{\partial x_1},
\cdots,\frac{\partial \vec f(x,u)}{\partial x_n})$.

\begin{prop} (\cite{Gory1}, Proposition 1.2 )
We have the equality
$$\mu =dim_{\bC} Q_f= dim \frac{{\mathcal O}_X}
{{\mathcal O}_X(f_1(x,u), \cdots, f_k(x,u))+ I_{C_0}}.$$
\label{prop1}
\end{prop}

Let us denote by $Cr(\vec F)$
the set of critical locus of the projection
$\pi: \bigcup_{t \in T} X_{t} \rightarrow S.$
That is to say
$$Cr(\vec F)=\{(x,u,t); (x,u) \in X_{t}, 
rank (  \frac{\partial \vec F(x,s)}{\partial x_1},
\cdots,\frac{\partial \vec F(x,s)}{\partial x_n} ) < k \}. \leqno(2.9)$$  We denote by
$D \subset S$ the image of projection $\pi(Cr(\vec F))$ which is
usually called discriminant set of the deformation $X_{t}$ of
projection. It is known that for the $R_+$-versal deformation, $D$
is defined by a principal ideal in ${\mathcal O}_S$ generated by a
single defining function $\Delta(s)$ \cite{Looi}. Under this
situation we define ${\mathcal O}_S-$ module of vector fields
tangent to the discriminant $D$ which is a sub-module of $Der_S$
the vector fields on $S$ with coefficients from ${\mathcal O}_S.$ 
\begin{dfn}
We define the {\rm logarithmic vector fields} associated to
$D$ as follows,
$$ Der_S(log\;D) =\{\vec v \in Der_S;\vec v(\Delta)\in
{\mathcal O}_S \cdot \Delta  \}.$$
We call that a meromorphic $p-$form
$\omega$ with a simple pole along $D$  belongs to the
 ${\mathcal O}_S$module of the {\rm logarithmic
differential forms } $\Omega^p_S(log\;D)$ associated to
$D$ iff the following two conditions are satisfied
$$ 1)   \Delta\cdot \omega \in \Omega^p_S,$$
$$ 2)   d\Delta\cdot \omega \in \Omega^{p+1}_S,$$
or equivalently
$$  \Delta\cdot d\omega \in \Omega^{p+1}_S.$$
\label{dfn1}
\end{dfn}

For the ${\mathcal O}_S$-module of the {\rm logarithmic
differential forms } the following fact is known.
\begin{thm}(See  \cite{Saito1} for the case $k=1$,
\cite{Looi}, \cite{Al} for the case $k$ general)
The module $ Der_S(log\;D)$  is a free ${\mathcal O}_S$-module
of rank $\mu$.  Furthermore there exists a $\mu$-tuple
of vectors $\vec v_1, \cdots, \vec v_\mu \in
Der_S(log\;D)$ such that
$$\Delta(s) = det (\vec v_1, \cdots, \vec v_\mu).$$
\label{thm2}
\end{thm}

\begin{prop}(see \cite{Terao} for the case $k=1$, \cite{Gory1}
for general $k$)

For every $\vec v_j \in Der_S(log\;D)$, $ 1 \leq j \leq \mu$,
there exists its lifting $\hat {\vec v}_j \in Der_{\tilde X \times
S}$ tangent to the critical set $Cr(\vec F).$ More precisely, the
following decomposition holds,
$$ \vec v_j (F_q(x,s))
= \sum_{p=1}^n h_{j,p}(x,s) \frac{\partial F_q}{\partial x_p}
+ \sum_{r=1}^k a_{jq}^{(r)}(x,s) F_r + b_{j,q}(x,s,\vec F),\; 1 \leq q \leq k$$
for some $h_{s,j}(x,s) \in {\mathcal O}_{\tilde X \times S}$,
 $b_{j,q}(x,s,\vec F) \in {\mathcal O}_{\tilde X \times S}
\otimes_{{\mathcal O}_{\tilde X \times S}} {m}^2_S. $
In this notation,
$$\hat {\vec v}_j = {\vec v}_j - \sum_{p=1}^n h_{j,p}(x,s)
\frac{\partial }{\partial x_p}.$$
Conversely, to every vector field $\hat {\vec v}_j \in Der_{\tilde X \times
S}$ tangent to the critical set  $Cr(\vec F)$ we can associate a vector field  $\vec v_j \in Der_S(log\;D)$ as its 
push down.
\label{prop2}
\end{prop}
This is a direct consequence of the preparation theorem (see \cite{Mal}).
Further on in this article we denote by $\vec v(F(x,s))$
the action of a vector field $\vec v \in \in Der_{\tilde X \times
S}$ on a function $F(x,s)$.

\begin{lem} (\cite{Gory1})
The discriminant $\Delta(s)$ defined in Theorem ~\ref{thm2}
can be expressed by a Weierstrass polynomial,
$$\Delta(s) = u^\mu + d_1(t) u^{\mu-1} + \cdots + d_\mu(t),$$
with $d_1(t) = \cdots =d_\mu(0)=0.$ 
\label{lem25}
\end{lem}
This can be deduced by another way by making use of
$(5.5)$ for the case of CI  $(5.1)$.
Namely we have $\Delta(s)= det P(s).$
From this lemma we deduce immediately the existence of an ``Euler''
vector field even for non-quasihomogeneous $\vec f(x,u)$
that plays essential r\^ole in the construction of the higher residue 
pairing by K.Saito\cite{Saito2}.

\begin{lem} (For $k=1$, see \cite{Saito2} (1.7.5))
There is a vector field 
$\vec v_1 = (u+\sigma_1^0(t)) \frac{\partial}{\partial u} + \sum^{\tau}_{i=1}
\sigma_1^i(t) \frac{\partial}{\partial s_i} \in  Der_S(log\;D)$
such that
$$ \vec v_1(\Delta(s)) = \mu \Delta(s).$$
\label{lem26}
\end{lem}
{\bf Proof}
It is clear that for a vector field $\vec v_1 \in 
Der_S(log\;D)$ with the component $(u+\sigma_1^0(t))
\frac{\partial}{\partial u}$ whose existence is guaranteed by Theorem
3,1 \cite{Gory1} ,
the expression 
$\vec v_1(\Delta(s))$ must be divisible by $\Delta(s)$.
In calculating the term of $\vec v_1(\Delta(s))$ 
that may contain the factor  $u^\mu $, we see that
$$\vec v_1(\Delta(s))= \mu u^\mu + {\tilde d}_1(t) u^{\mu-1} + \cdots + {\tilde d_\mu}(t).$$
Thus we conclude that ${\tilde d}_i(t) = \mu  d_i(t),$ $1 \leq i \leq \mu.$ 
{\bf Q.E.D.}

Now we introduce the filtered ${\mathcal O}_{S}$-module
of fibre integrals ${\mathcal H}^{(\vec \lambda)}$
for a multi-index $\vec \lambda =(\lambda_1, \cdots,
\lambda_k)\in (\bZ_{<0})^k$.
$$ I_\phi^{\vec \lambda}(s) = \int_{t(\gamma)} \phi(x,s)
F_1(x,s)^{\lambda_1} \cdots
F_k(x,s)^{\lambda_k} dx,$$
for $\phi(x,s)  \in {\mathcal O}_{\tilde X \times S}.$
Let us denote by $X^{(q)} := \{x \in \tilde X; F_q(x,s)=0 \} $
a smooth hypersurface defined for $s \not \in D.$
In this situation 
we define the Leray's tube operation isomorphism (see \cite{Vas}, \cite{HT}),
$$ \begin{array}{cccc}
t:& H_{n-k}(\cap_{q=1}^k X^{(q)})& \rightarrow & H_{n}(\tilde X \setminus
\cup_{q=1}^k X^{(q)}),\\ 
&\gamma & \mapsto & t(\gamma).
\end{array}$$
The concrete construction of the operation $t$ can be described as follows.
First we consider the coboundary isomorphism of the homology groups,
$$ \delta: H_{n-k}(\cap_{q=1}^k X^{(q)}) \rightarrow  H_{n-k+1}(\cap_{q=2}^k X^{(q)} \setminus X^{(1)}).$$
A cycle $\gamma$ in $\cap_{q=1}^k X^{(q)}$ is mapped onto a cycle $\delta(\gamma)$ of one higher dimension
that is obtained as a $S^1$ bundle over $\gamma$.
Repeated application of  $ \delta$ yields an interated coboundary homomorphism,
$$ H_{n-k}(\cap_{q=1}^k X^{(q)}) {\rightarrow^\delta} H_{n-k+1}(\cap_{q=2}^k X^{(q)} \setminus X^{(1)}) 
{\rightarrow^\delta}
\cdots$$
$$ \cdots  {\rightarrow^\delta}H_{n-1}( X^{(k)} \setminus \cup_{q=1}^{k-1} X^{(q)}) 
{\rightarrow}^{\delta}H_{n}( \tilde X \setminus \cup_{q=1}^{k} X^{(q)}).$$
The Leray's tube operation is a $k-$time iterated $\delta$ homomorphism i.e.
$t = \delta^m$. The Froissart decomposition theorem (\cite{HT}, \S 6-3)
shows that the collection of all cycles of $H_{n}( \tilde X \setminus \cup_{q=1}^{k} X^{(q)})$
are obtained by the application of iterated $\delta$ homomorphism operations to
the cycles from $H_{n-p}( \tilde X \cap
X^{(q_1)}\cap X^{(q_2)} \cdots \cap X^{(q_p)})$, $p=0, \cdots, k$.

Let us denote by
$\Phi$ the $\bC$ vector space  $\frac{{\mathcal O}_{X_0}}{I_{C_0}}$
whose $\bC-$ dimension is equal to $\mu$ after the Proposition
~\ref{prop1}. We denote its basis by $(\phi_0(x,u), \cdots,
\phi_\tau(x,u))$

Now let us introduce a notation of the multi-index $ -\b1 =(-1,
\cdots-1) \in (\bZ_{<0})^k$. We consider a vector of fibre
integrals $\bI_{\Phi}:= ^t(I_{\phi_0}^{(-\b1)}(s), \cdots,
I_{\phi_\tau}^{(-\b1)}(s)).$ The following theorem for $k=1$
has been anounced in \cite{Saito3} (4.14) without proof.

\begin{thm}

1.For every $\vec v \in Der_S(log\;D)$,
we have the following inclusion relation
$$\vec v :{\mathcal H}^{(-\b1)} \hookrightarrow {\mathcal H}^{(-\b1)}.$$
That is to say for every  $\vec v_j \in Der_S(log\;D)$, there
exists a $\mu \times \mu$ matrix with holomorphic entries
$B_j(s)\in  End(\bC^\mu) \otimes {\mathcal O}_{S}$ such that
$$ \vec v_j(\bI_{\Phi}) = B_j(s) \bI_{\Phi}, 1 \leq j \leq \mu.$$

2. The vector of fibre integrals $\bI_{\Phi}$ satisfies the following
Pfaff system of Fuchsian type
$$ d\bI_{\Phi} = \Omega \cdot \bI_{\Phi},$$
 for some
$\Omega \in End(\bC^\mu) \otimes_{{\mathcal O}_S} \Omega^1_S (log\;D).$

\label{thm3}
\end{thm}

{\bf Proof} As for the proof of 1, we remark the following
equality that yields from Proposition ~\ref{prop2},
$$\vec v_j \left( \int_{t(\gamma)} \phi(x,u)
F_1(x,s)^{-1} \cdots F_k(x,s)^{-1} dx \right) = \int_{t(\gamma)} \vec
F^{-\b1} d(\phi(x,u) \sum_{p=1}^n (-1)^{p-1} h_{j,p}(x,s) dx_1
\pvup dx_n) + $$ $$+\int_{t(\gamma)} \vec F^{-\b1} \phi(x,u)
(\sum_{q=1}^k \sum_{r=1}^k a_{j,q}^{(r)} F_r F_q^{-1})dx
+\int_{t(\gamma)} \vec F^{-\b1}  \sum_{q=1}^k
F_q^{-1}b_{j,q}(x,u,t,\vec F)dx $$
$$=\int_{t(\gamma)} \vec
F^{-\b1} d(\phi(x,u) \sum_{p=1}^n (-1)^{p-1} h_{j,p}(x,s) dx_1
\pvup dx_n)+\int_{t(\gamma)} \vec F^{-\b1} \phi(x,u) (
\sum_{r=1}^k a_{j,r}^{(r)}(x,s))dx,$$ which evidently belongs
to ${\mathcal H}^{(-\b1)}$. The last equality can be explained by
the vanishing of the integral
$$\int_{t(\gamma)} F_1^{-1} \cdots 1 \rvup   F_q^{-2}\qvup F_k^{-1} \phi(x,u)
( a_{j,q}^r  )dx =0,$$ because of the lack of the residue along
$F_r(x,s)=0$ and
$$\int_{t(\gamma)} \vec F^{-\b1} F_{q_1} F_{q_2} F_q^{-1} \phi(x,u)
( b^{0}_{j,q}(x,s)  )dx =0,$$ in view of the lack of at least one
of residues either along $F_{q_1}=0$ or along $ F_{q_2}=0$.
These equalities are derived from the property of the Leray's tube
$t(\gamma)$ which needs codimension $k$ residue to give rise to a non-zero
integral.

2.
Let us rewrite the relations obtained in 1. into the form,
$$dI_{\phi_q}^{(-\b1)}=
\sum_{r=1}^\mu \omega_{q,r} I_{\phi_r}^{(-\b1)},$$
for some $\omega_{q,r} \in \Omega_S^1(-D)$ meromorphic 1-forms with poles
along $D.$ These $\omega_{q,r}$ satisfy the following relations,
$$\vec v_j(I_{\phi_q}^{(-\b1)})=
\langle \vec v_j, dI_{\phi_q}^{(-\b1)} \rangle =\langle \vec v_j,
\sum_{r=1}^\mu \omega_{q,r} I_{\phi_r}^{(-\b1)} \rangle\;\; 1 \leq j,q \leq 
\mu.$$ If
$\langle \vec v_j, \omega_{q,r}\rangle \in {\mathcal O}_S$  for
all $\vec v_j \in Der_S(log\;D)$ $1 \leq j \leq \mu$ then
$\omega_{q,r} \in \Omega^1_S(log\;D)$ in view of the
Theorem~\ref{thm2}. {\bf Q.E.D.}

Let us introduce a filtration as follows
${\mathcal H}^{(\lambda)}= \bigoplus_{\lambda_1 + \cdots + \lambda_k = \lambda}
{\mathcal H}^{(\vec \lambda)}$.
For this rough filtration we have the following generalisation of the Griffiths' transversality
theorem (\cite{GSch} Theorem 3.1).
\begin{cor}
For every $\vec v \in Der_S(log\;D)$,
we have the following inclusion relation
$$\vec v : {\mathcal H}^{(\lambda)}\hookrightarrow 
{\mathcal H}^{(\lambda)}.$$
\end{cor}
{\bf Proof}
For $\partial_{s_j} I_{\Phi} \in {\mathcal H}^{(-k-1)}$
and $\vec v_\ell \in Der_S(log\;D)$ we have
$$\vec v_\ell(\partial_{s_j} I_{\Phi})
= [\vec v_\ell,\partial_{s_j}]I_{\Phi} + \partial_{s_j} \vec v_\ell(I_{\Phi})$$
$$=[\vec v_\ell,\partial_{s_j}]I_{\Phi}+ \partial_{s_j} (B_\ell(s) I_{\Phi})
=[\vec v_\ell,\partial_{s_j}]I_{\Phi}+ (\partial_{s_j} B_\ell(s)) I_{\Phi}
+B_\ell(s) (\partial_{s_j}I_{\Phi}).$$
As the commutator
$[\vec v_\ell,\partial_{s_j}]$ is a first order operator, the term above
$[\vec v_\ell,\partial_{s_j}]I_{\Phi}$ belongs to  ${\mathcal H}^{(-k-1)}.$
The term $\partial_{s_j} B_\ell(s) I_{\Phi} \in {\mathcal H}^{(-k)}$
again belongs to  ${\mathcal H}^{(-k-1)}.$ Thus we see
$\vec v_\ell(\partial_{s_j} I_{\Phi}) \in {\mathcal H}^{(-k-1)}$.
In an inductive way, for any $\lambda \leq -k$ we prove the statement.

{\bf Q.E.D.}

{ \center{\section{Multiplication table and the logarithmic vector fields}} }

We consider a miniversal deformation of  a
mapping $\vec f^{(0)}(x)$
which can be written down in the following special form
for $s=(u,t) \in S$,
$$\vec F(x,s) = \vec f^{(0)}(x) + \sum_{\ell=1}^\tau t_\ell \vec e_\ell(x)
+ u\vec e_0(x) =\left (\begin {array}{c}
F_1(x,t)-u\\
F_2(x,t)\\
\vdots\\
F_k(x,t)\\
\end {array}\right ),\leqno(3.1)$$
for
$$\{ \vec e_0(x), \cdots , \vec e_\tau(x)\} \in Q_f,$$
where $\vec e_0(x) = \;^t(-1, 0, \cdots,0).$
One may consult \cite{Looi} (6.7) to see that $\vec F(x,s)$
really gives a miniversal deformation of  $\vec f^{(0)}(x) $
by virtue of the definitions $(2.3)$, $(2.5).$
Let us fix a basis
$\{ \phi_0(x), \cdots, \phi_\tau(x) \}$
of the space  $\Phi:= \frac{{\mathcal O}_{X}}
{I_{C_0}+{\mathcal O}_X(f_1(x)-u, f_2(x), \cdots, f_k(x))}.$
We remark here that the basis of $\Phi$ can be represented by
functions from ${\mathcal O}_{\tilde X}$ as we can erase the variable $u$ by the relation
$f_1(x)=u$ in $\Phi$.
It turns out that we can regard
$\{ \phi_0(x), \cdots, \phi_\tau(x) \}$ as a free basis of the
${\mathcal O}_S$ module $\Phi(s)$ treated in
the Proposition ~\ref{prop50} below.
Under these circumstances,
we introduce holomorphic functions $\tau^\ell_{i,j}(s) \in {\mathcal O}_S$
in the following way.
$$ \phi_i(x) \vec e_j(x) \equiv
\sum_{\ell=0}^\tau \tau^\ell_{i,j}(s) \vec e_\ell(x)
\;\;mod ({\mathcal O}_{\tilde X \times S}
\langle \frac{\partial \vec F(x,s)}{\partial x_1}, \cdots,
\frac{\partial \vec F(x,s)}{\partial x_n}\rangle).\leqno(3.2) $$
The functions $\tau^\ell_{i,j}(s) \in {\mathcal O}_S$ 
exist due to the versality of the deformation
$\vec F(x,s).$
We denote by
$$T_j(s) =\left( \tau^\ell_{i,j}(s)   
\right)_{0 \leq j,\ell \leq \tau}, \leqno(3.3) $$
a $\mu \times \mu$ matrix which is called 
the matrix of {\bf multiplication table}.
We denote the discriminant associated to this deformation
by $D  \subset S$.

Further on we will make use of the abbreviation $mod (d_x \vec F(x,s))$
instead of making use of the expression $mod 
({\mathcal O}_{\tilde X \times S}\langle 
\frac{\partial \vec F(x,s)}{\partial x_1}, \cdots,
\frac{\partial \vec F(x,s)}{\partial x_n}\rangle).$

After Proposition ~\ref{prop2} the vector field $\vec v_1$
constructed in Lemma ~\ref{lem26} has its lifting
$\hat {\vec v_1} \in Der_{\tilde X \times S}.$ Let us denote
by $\check {\vec v_1}= \hat{\vec v_1} - \vec v_1$ $\in 
{\mathcal O}_{\tilde X \times S} \otimes Der_{\tilde X}$.
$$  \hat{\vec v_1}(\vec F(x,s))\cdot \phi_i(x)
= \check {\vec v_1}(\vec f^{(0)}(x))\cdot \phi_i(x) + \sum_{\ell=0}^\tau
\vec v_1(s_\ell) \vec e_\ell(x)\phi_i(x) + \sum_{\ell=0}^\tau
s_\ell (\check{\vec v_1}e_\ell(x))\phi_i(x) $$
$$ \equiv \sum_{\ell=0}^\tau
{\vec v_1}(s_\ell) \vec e_\ell(x)\phi_i(x) \;\;\;mod (d_x \vec F(x,s)).$$
\begin{lem}
There exists a vector valued function
$M(x,\vec F(x,s))$
$\in$ $({\mathcal O}_{{\tilde X}\times \bC^k})^k$
such that $$ \hat{\vec v}_1(\vec F(x,s))
\equiv M (x,\vec F(x,s)) \;\;\;mod (d_x \vec F(x,s)),$$ with
$$M(x,\vec F(x,s)) = M^0 \cdot \vec F(x,s) +M^1(x,\vec F(x,s)),$$
where  $ M^0 \in GL(k, \bC)$: a non-degenerate matrix and $M^1(x,\vec F(x,s))$
$\in ({\mathcal O}_{\tilde X}\otimes m^2_S)^k$.
Especially the first row of $ M^0=$ 
$(1,0,\cdots,0).$
\label{lem31}
\end{lem}

{\bf Proof}
First of all we remember a theorem due to \cite{Gr} \S 1.1, 
\cite{Saito3} Proposititon 2.3.2 which states that the Krull 
dimension of the ring of 
holomorphic functions
on the critical set
$Cr(\vec F)$ is equal to $\mu-1$ and this ring is a Cohen- Macaulay ring.
Let us denote by $L=\;_nC_k$.  We  have  $(k+L)$tuple of $k \times k-$ minors 
$j_{k+1}(x,s)$ $\cdots$ $j_{k+L}(x,s)$ of the matrix
$( \frac{\partial}{\partial x_1} 
\vec F(x,s), \cdots,\frac{\partial}{\partial x_n} \vec F(x,s)) $
such that
$$Cr(\vec F)= V(\langle F_1(x,s), \cdots, F_k(x,s),j_{k+1}(x,s), \cdots, j_{k+L}(x,s)\rangle ).$$
The lemma 
~\ref{lem26} yields that the lifting $\hat {\vec v}_1$
of the vector field $\vec v_1$ 
satisfies the relations,
$$ 
\langle F_1(x,s), \cdots, F_k(x,s),j_{k+1}(x,s), \cdots, j_{k+L}(x,s)\rangle $$
$$ = 
\langle\hat{\vec v}_1(F_1(x,s)), \cdots, \hat{\vec v}_1(F_k(x,s)),
\hat{\vec v}_1(j_{k+1}(x,s)), \cdots, \hat{\vec v}_1(j_{k+L}(x,s))\rangle.$$
As it has been seen from the above Proposition~\ref{prop2}, the vector $\hat{\vec v}_1$
is tangent to $Cr(\vec F)$.  If the above equality does not hold,
it would entail the relation
$$ D =\{s \in S; \Delta(s)=0\} {\subsetneqq} \pi (V(\langle\hat{\vec v}_1(F_1(x,s)), \cdots, \hat{\vec v}_1(F_k(x,s)),
\hat{\vec v}_1(j_{k+1}(x,s)), \cdots, \hat{\vec v}_1(j_{k+L}(x,s))\rangle)),$$
after elimination theoretical consideration.
This yields 
$$ \hat {\vec v}_1(F_q(x,s)) = \sum_{\ell=1}^k C_q^\ell F_\ell(x,s)+m_q(x,\vec F) 
+ \sum_{\ell=k+1}^{k+L} C_q^\ell j_\ell(x,s), 1 \leq q \leq k,  $$
$$ \hat {\vec v}_1(j_p(x,s)) = 
 \sum_{\ell=k+1}^{k+L} C_p^\ell j_\ell(x,s) + m_p(x,\vec F), k+1 \leq p \leq k+L,$$
for $m_r(x,\vec F) \in {\mathcal O}_{\tilde X}\otimes m^2_{S}$, $1 \leq r 
\leq k +L $
and some constants $C_q^\ell, 1 \leq \ell \leq k.$
First we see that the expression $ \hat {\vec v}_1(j_p(x,s))$
cannot contain terms of $F_q(x,s)$ like $F_q(0,s)$ in view of the situation that
the versality of the deformation makes all linear in $x$ variable terms dependent
on some of deformation parameters.
Secondly the non-degeneracy of the matrix $M^0:= (C_q^\ell)_{ 1 \leq q, \ell \leq k}$
is necessary so that the above equality among ideals holds. 

From this relation and the preparation theorem, we see
$$\hat {\vec v}_1(\vec F(x,s))= M^0 \cdot \vec F(x,s) +M^1(x,\vec F(x,s))+ h_{1,1}(x,s)
\frac{\partial \vec F(x,s)}{\partial x_1}+ \cdots + h_{1,n}(x,s)
\frac{\partial \vec F(x,s)}{\partial x_n},$$
  with $M^1(x,\vec F(x,s))= \;^t(m_1(x,\vec F), \cdots, m_k(x,\vec F)) \in 
({\mathcal O}_{\tilde X}\otimes m^2_{S})^k$.

More precisely we can state that $C_1^1=1,$
$C_1^\ell =0, 2 \leq \ell \leq k$. 
The dependence of some cofficients of $\hat {\vec v}_1$ on $F_i(x,t)$
is necessary so that $C_1^\ell \not= 0$ for some $2 \leq \ell \leq k.$
But this is impossible because
if not it would mean that some of the coefficients of 
$\hat {\vec v}_1$ contains factor $F_2(x,s), \cdots, F_k(x,s)$
that contradicts the construction of $\hat {\vec v}_1$ in Proposition
~\ref{prop2}.
This can be seen from the fact that the expressions
$\frac{\partial F_1(x,s)}{\partial x_1},$ $\cdots,$
$\frac{\partial  F_1(x,s)}{\partial x_n},$
$\frac{\partial F_1(x,s)}{\partial s_1},$ $\cdots,$
$\frac{\partial F_1(x,s)}{\partial s_\mu}$
do not contain the deformation parameters
present in the polynomials $F_2(x,s), \cdots, F_k(x,s).$
{\bf Q.E.D.}

\begin{lem}
A basis
of  logarithmic vector fields
$\vec v_0, \cdots, \vec v_\tau  \in Der_S(log\;D)$
can be produced from the functions $\sigma^\ell_{i}(s)$ defined as follows,
$$\hat {\vec v}_1(\vec F(x,s))\cdot \phi_i(x) = M(x,\vec F(x,s)) \cdot \phi_i(x) = \sum_{\ell=0}^\tau 
\sigma^\ell_{i}(s) {\vec e_\ell} + \check{\vec v_i}(\vec F(x,s))
$$
$$ \equiv \sum_{\ell=0}^\tau \sigma^\ell_{i}(s) {\vec e_\ell}  \;\;
mod (d_x \vec F(x,s)),$$
where the vector valued fucntion $M(x,\vec F(x,s))$ 
denotes the one defined in the Lemma ~\ref{lem31}
and $\check{\vec v_j}=\sum_{p=1}^n
h_{j,p}(x,s)\frac{\partial  }{\partial x_p}$
is a certain vector field with holomorphic coefficients.
 
\label{lem32}
\end{lem}

{\bf Proof}

We remark the following relation,
$$\hat {\vec v}_1(\vec F(x,s)) \phi_i(x) =
\check {\vec v}_1(\vec f^{(0)}(x))\phi_i(x)
+ \sum_{j=1}^\mu \vec v_1(s_j) \vec e_j(x)\phi_i(x)+\sum_{j=1}^\mu s_j \check {\vec v}_1
(\vec e_j(x))\phi_i(x)$$
  $$ \equiv \sum_{j=0}^\tau \vec v_1(s_j) \vec e_j(x)\phi_i(x)\;\;
mod (d_x \vec F(x,s)).$$
The relation $(3.2)$ above entails,
$$M(x,\vec F(x,s)) \cdot \phi_i(x)
\equiv \sum_{\ell=0}^\tau \sum_{j=0}^\tau \vec v_1(s_j)
\tau_{i,j}^\ell(s)  \vec e_\ell(x)\;\; mod (d_x \vec F(x,s)).$$
As $\phi_i(x)$ can be considered to be  
a basis of ${\mathcal O}_S$ module  $\Phi(s)$ above
(see Proposition~\ref{prop50}),
vectors $(\sigma^0_{i}(s), \cdots, \sigma^\tau_{i}(s)),$ 
$0 \leq i \leq \tau$ are ${\mathcal O}_S$ linearly 
independent at each generic point $S \setminus D$.
If we put
$$\sigma^\ell_{i}(s)=\sum_{j=0}^\tau \vec v_1(s_j)
\tau_{i,j}^\ell(s),$$ 
then the vector field $\hat{\vec v_i}\in Der_{\tilde X \times S}$
$$ \hat{\vec v_i}= \sum^\tau_{\ell=0} \sigma^\ell_{i}(s)
\frac{\partial  }{\partial s_\ell} + \phi_i(x) \check {\vec v}_1,$$
is tangent to $Cr(\vec F)$. The only non-trivial relations
that may arise between $ \check{\vec v_i}$ and $ \check{\vec v_{i'}}$ $i \not = i'$
is 
$$ \phi_i(x) \check{\vec v_{i'}}= \phi_{i'}(x) \check{\vec v_i}.$$
These vectors give rise to the same push down vector field in
 $Der_S(log\;D)$.
Namely,
$$ \pi_\ast(\phi_i(x) \hat{\vec v_{i'}})= \pi_\ast(\phi_{i'}(x) \hat{\vec v_{i}})=
\sum^\tau_{j=0}\sum^\tau_{\ell=0} R^\ell_{i,i',j}(s)\frac{\partial  }{\partial s_\ell},$$
for the coefficients $R^\ell_{i,i',j}(s)$ determined by
$$ \sum^\tau_{j=0} \vec v_1(s_j)\phi_i(x) \phi_{i'}(x)\vec e_j(x) \equiv
\sum^\tau_{j=0}\sum^\tau_{\ell=0} R^\ell_{i,i',j}(s)\vec e_\ell(x) mod(d_x \vec F(x,s)).$$
This means that $\hat{\vec v_0}, \cdots, \hat{\vec v_{\tau}}$ form a free basis of
$Der_{\tilde X \times S}(Cr(\vec F))$ hence  ${\vec v_0}, \cdots, {\vec v_{\tau}}$
that of  $ Der_S(log\;D).$
{\bf Q.E.D.}

This lemma gives us a correspondence between  $ \phi_i(x)\in \Phi$
and ${\vec v_i} \in Der_S(log\;D),$ therefore 
it is quite natural to expect that the 
mixed Hodge structure on $\Phi$ would induce that on 
$Der_S(log\;D)$, and would hence contribute to describe $B_i(s)$
of Theorem ~\ref{thm3}, 1 in a precise manner. 
A good understanding of this situation is indispensable to characterize 
the rational monodromy of solutions to the Gauss-Manin system in terms of the 
mixed Hodge structure on $\Phi$. Confer to  Proposition ~\ref{prop53} below.
 
We formulate the lemma~\ref{lem32} into the following form
(see \cite{bru} Theorems A2, A4, \cite{Saito1} (3.19), \cite{Saito3} (4.5.3)
Corollary 2  for $k=1$ and \cite{Looi} (6.13), \cite{Gory1} Theorem 
3.2 for $k$ general).

\begin{prop}
There exist holomorphic functions $w_j(s) \in {\mathcal O}_S,$ $0 \leq j \leq \tau$
such that the components of the matrix
$$ \Sigma(s):= \sum_{j=0}^\tau  w_j(s) T_j(s), \leqno(3.4) $$
give rise to  a basis of  logarithmic vector fields
$\vec v_0, \cdots, \vec v_\tau  \in Der_S(log\;D)$.
Namely, if we write $\Sigma(s)=\left( \sigma^\ell_{i}(s)   \right)_{0 
\leq i,\ell \leq
\tau},$
then the expression
$$\vec v_i= \sum_{\ell =0}^\tau 
\sigma^\ell_{i}(s)\frac{\partial}{\partial s_\ell},
\leqno(3.5) $$
consists a base element of the ${\mathcal O}_S$
module $ Der_S(log\;D).$ 
\label{prop31}
\end{prop}

Especially in the case of quasihomogeneous singularity
$\vec f(x,u)$ we have the following simple description of the vector field
that can be deduced from Lemma ~\ref{lem32}.
To do this, it is enough to remark that the vector field $\vec v_1$
is the Euler vector field by definition and
$\vec v_1(s_r) = \frac{w(s_r)}{w(s_0)}s_respectr,$
where $w(s_j)$ denotes the quasihomoeneous weight  of the variable $s_j.$
\begin{prop}(\cite{Gory2} Theorem 2.4)
In the case of quasihomogeneous singularity
$(2.1)$,
the basis $(3.5)$ of $Der_S(log\;D)$ can be calculated by
$$\sigma_{i}^{\ell}(s) = \sum_{j=0}^\tau w(s_j) s_j \tau^\ell_{i,j}(s).$$
Furthermore, the vector valued function  $M(x,\vec F(x,s))$
of Lemma ~\ref{lem31} has the expression,
$$ M(x,\vec F(x,s))=M^0 \cdot \vec F(x,s)= diag\;(w(f_1), \cdots, w(f_k))\cdot
\vec F(x,s).$$
\label{prop34}
\end{prop}


{ \center{\section{Multiplication table and the topology of 
real hypersurfaces}} }

In this section we continue to consider the situation where $\mu=\tau+1$
for $k=1$ in $(2.5).$
We associate to  the versal deformation of the hypersurface singularity
$$F(x,s) = f(x) + \sum_{i=0}^\tau s_i e_i(x),\leqno(4.1)$$
 the following matrix
$\Sigma(s)=(\sigma_i^\ell(s))_{0 \leq i,\ell \leq \tau}$ 
after the model $(3.2)$,

$$ F(x,s) e_i(x) = \sum_{\ell=0}^\tau \sigma^\ell_{i}(s) e_\ell(x)
\;\;mod (d_x F(x,s)). \leqno(4.2)$$

$$ e_i(x) e_j(x) \equiv \sum_{\ell=0}^\tau \tau^\ell_{i,j}(t)  e_\ell(x)
\;\;mod (d_x F(x,s)).\leqno(4.3)$$
Further on we make use of the convention
$e_0(x)=1$ and $s=(s_0,t).$ We denote the deformation
parameter space $t \in T= (\bC^\tau,0)$.

We recall the Milnor ring for $k=1$ whose analogy has been
introduced in $(2.5)$ (and in the case
$k$ general, $\Phi(s)$
will be introduced in Proposition ~\ref{prop50}), 
$$ Q_F:= \frac{{\mathcal O}_{\tilde X \times S}}
{{\mathcal O}_{\tilde X \times S}\langle 
\frac{\partial F(x,s)}{\partial x_1}, \cdots,
\frac{\partial F(x,s)}{\partial x_n}\rangle}.$$
We introduce the Bezoutian matrix $B^F(s)$
whose $(i,j)$ element is defined by the trace of
the multiplication action $F(x,s)e_i(x) e_j(x)\cdot$ 
on the Milnor ring $Q_F,$
$$F(x,s)e_i(x) e_j(x)  \equiv (\sum_{c=0}^\tau 
\sigma_i^c(s) e_c(x))e_j(x) $$
$$ \equiv \sum_{c=0}^\tau\sigma_i^c(s)(\sum_{r=0}^\tau 
\tau^r_{c,j}(t) e_r(x))
mod (d_x F(x,s)).$$
For the sake of simplicity we will use the following notation,
$$\tau^r(t) = ( \tau^r_{c,b}(t) )_{0 \leq c,b \leq \tau}. \leqno(4.4)$$
To clarify the structure of the Bezoutian matrix $B^F(s)$
we introduce a matrix
$$T(t) = \left( \sum_{r=0}^\tau \zeta_r(t)  
\tau^r(t) \right ), \leqno(4.5)$$
with the notation
$$ \zeta_r(t) = tr(e_r(x) \cdot)= 
\sum_{\ell=0}^\tau \tau^\ell_{r,\ell}(t). \leqno(4.6)$$
The $(i,j)$ element of the matrix $T(t)$ $(4.5)$ equals to
$tr(e_i(x) e_j(x) \cdot)$ on the Milnor ring $ Q_F$.
It is possible to show  that $\{t \in T;det(T(t))=0\}$
coincides with the bifurcation set of $F(x,s)$
outside the Maxwell set (see Proposition ~\ref{prop52} below).
Thus we get the Bezoutian matrix
$$B^F(s) = \Sigma(s) \cdot T(t). \leqno(4.7)$$
Following statement is a simple application
of Morse theory to the multiplication table
see  \cite{Sz1} Theorem 2.1. From here on we assume that
$|s|$ is small enough and  denote by $\tilde X
=\{x \in \bC^n; |x| \leq \delta\}$ a closed ball 
such that all critical points of 
$F(x,s)$ are located inside $\tilde X.$

\begin{prop}
$sign\; \Sigma(s)\cdot T(t)$ =$\{$ number of real critical points 
with respect to the variables $x$ in 
$F(x,s)>0$,  $x \in \tilde X \cap \bR^n\}$
-$\{$ number of real critical points with respect to the variables $x$ 
in $F(x,s)<0$, $x \in \tilde X \cap \bR^n
\}$.
Here $sign(A)$ denotes the signature of a symmetric 
matrix $A$ i.e. the difference between the number of positive and
negative eigenvalues.
\label{prop41}
\end{prop}

Let us denote by $h(x,t)$
the determinant of the Hessian
$$ h(x,t) := det \langle \frac{ \partial^2 F(x,s)}{ \partial x_i \partial x_j}
\rangle_{1 \leq i,j \leq n}.$$
We associate the following $\mu$ holomorphic functions $h_0(t), \cdots, h_\tau(t) 
\in {\mathcal O}_{S}$
to the function $h(x,t)$,
$$h(x,t)\equiv \sum_{\ell=0}^\tau h_\ell(t) e_\ell(x)
\;\;mod (d_x F(x,s)).\leqno(4.8)$$
Further by means of $(4.7)$ we introduce the matrix
$$ B^H(t) :=\sum_{\ell=0}^\tau \eta^\ell(t)  \tau^\ell(t), \leqno(4.9)$$
where
$$
\left (\begin {array}{c}
\eta^0(t)\\
\vdots\\
\eta^\tau(t)\\
\end {array}\right )=
T(t)\cdot
\left (\begin {array}{c}
h_0(t)\\
\vdots\\
h_\tau(t)\\
\end {array}\right ).$$
We consider the matrix $B^{HF}(s)=(\cdot)_{0 \leq a,b \leq \tau}$ 
whose $(a,b)-$element is defined by the trace of the following expression
on the Milnor ring $Q_F$,
$$h(x,t) F(x,s) e_a(x) e_b(x)
\equiv (\sum_{\ell=0}^\tau h_\ell(t) e_\ell(x))(\sum_{c=0}^\tau 
\sigma_a^c(s)\sum_{m=0}^\tau \tau^m_{c,b}(t)e_m(x)) 
\leqno(4.10)$$
$$\equiv \sum_{\ell=0}^\tau  \sum_{c=0}^\tau \sum_{m=0}^\tau 
h_\ell(t) \sigma_a^c(s)\tau^m_{c,b}(t) e_\ell(x) e_m(x)$$
$$\equiv \sum_{\ell=0}^\tau  \sum_{c=0}^\tau \sum_{m=0}^\tau 
h_\ell(t) \sigma_a^c(s)\tau^m_{c,b}(t) \sum_{r=0}^\tau \tau^r_{\ell,m}(t)
e_r(x)\;\;\;mod( d_xF(x,s)).$$
If we take the trace of this, we get
$$  \sum_{c=0}^\tau \sigma_a^c(s) 
\sum_{m=0}^\tau \sum_{\ell=0}^\tau h_\ell(s) \tau^m_{c,b}(t)
(\sum_{r=0}^\tau \tau^r_{\ell,m}(t) \zeta_r(t)).$$
After $(4.8)$ and $(4.9)$ this matrix has the
following expression,
$$ B^{HF}(s) = \Sigma(s)\cdot B^H(t). \leqno(4.11)$$

We consider the following closures of semi-algebraic sets,
$$W_{= 0} := \{x \in \tilde X \cap \bR^n; F(x,s) = 0  \},$$
$$W_{\geq 0} := \{x \in \tilde X \cap \bR^n; F(x,s) \geq 0  \},
W_{\leq 0} := \{x \in \tilde X \cap \bR^n; F(x,s) \leq 0  \}.$$
\begin{thm}
The following expression of the Euler characteristics for $W_{\ast}$
holds,
$$\chi(W_{\geq 0}) - \chi(W_{=0}) = \frac{sign( B^{H}(t))+
sign ( B^{HF}(s))}{2}.$$
$$\chi(W_{\leq 0}) - \chi(W_{=0}) = (-1)^n\frac{sign( B^{H}(t))-
sign ( B^{HF}(s))}{2}.$$
\label{thm42}
\end{thm}
{\bf Proof}

After Szafraniec \cite{Sz1}, or simply applying Morse theory
to the real fibres of $F(x,s)$, we have the following equalities,
$$ \sum_{x \in {\rm critical\; points \;of} F(x,s) } (sgn \;h(x,t))$$
$$=sign\langle tr(h(x,t)e_i(x)\cdot e_j(x) \cdot)\rangle_{1 \leq i,j \leq n}
=\sum_{x \in {\rm critical\; points \;of} F(x,s) } (-1)^{\lambda(x)}.$$
Here we denoted by $tr(h(x,t)e_i(x)\cdot e_j(x)\cdot)$
the trace of a matrix defined by the multiplication 
by  the element $h(x,t)e_i(x)\cdot e_j(x)$ 
considered $mod( d_xF(x,s))$ for the basis $e_i(x), 1 \leq i \leq \mu.$
$$ \sum_{x \in {\rm critical\; points \;of} F(x,s) } (sgn \;h(x,t)) (sgn \;F(x,s))$$
$$=sign\langle tr(h(x,t)F(x,s) e_i(x)\cdot e_j(x)\cdot)
\rangle_{1 \leq i,j \leq n}
=\sum_{x \in {\rm critical\; points \;of} F(x,s) } (-1)^{\lambda(x)} (sgn \;F(x,s)).$$
We denoted by $tr(h(x,t)F(x,s)e_i(x)\cdot e_j(x)\cdot)$
the trace of a matrix defined by the multiplication 
by  the element $h(x,t)F(x,s)e_i(x)\cdot e_j(x)$ 
considered $mod( d_xF(x,s))$ for the basis $e_i(x), 1 \leq i \leq \mu.$
The exponent $\lambda(x)$ is the Morse index of the function $F(x,s)$ at $x$
and $sgn \;h(x,t)=(-1)^{\lambda(x)}$. {\bf Q.E.D.}

{ \center{\section{Topology of real complete intersections}} }

Let us reconsider the situation $(3.1)$
for the deformation of the CI,
$$
\vec F(x,u,t) =
\left (\begin {array}{c}
F_1(x,t)-u\\
F_2(x,t)\\
\vdots\\
F_k(x,t)\\
\end {array}\right ), \leqno(5.1) $$
with $s = (u,t ) \in S.$
Define the ideal
$I_{C_0}(t) \subset {\mathcal O}_{{\tilde X} \times S}$
generated by $k \times k $ minors of the marix
$ (\frac{\partial \vec F(x,0,t)}{\partial x_1},
\cdots,\frac{\partial \vec F(x,0,t)}{\partial x_n})$.

We have the following isomorphisms
$$ \Phi =
\frac{ {\mathcal O}_{{X}}}{{\mathcal O}_{{ X}}
\langle f_1(x)-u,f_2(x),\cdots, f_k(x) \rangle  + I_{C_0}(0)} $$
$$ \cong
\frac{ {\mathcal O}_{{\tilde X}}}{{\mathcal O}_{\tilde X}
\langle f_2(x),\cdots, f_k(x) \rangle + I_{C_0}(0)}, \leqno(5.2)$$
where $I_{C_0}(0)$ is the corresponding ideal in 
${\mathcal O}_{\tilde X}$.
The dimension of this space is equal to $\mu$ introduced 
in Proposition ~\ref{prop1}.
As for this number we remember that it can be expressed by means of
the Milnor number of the singularity $X_1:=\{x\in \tilde X; f_2(x)=\cdots=f_k(x)=0  \}$
and the Milnor number of the function $f_1$ restricted on $X_1$ i.e.
that of the singularity $\tilde X_0:=\{x\in \tilde X; f_1(x)=f_2(x)=\cdots=f_k(x)=0  \}$,
$$ \mu = \mu(X_1) + \mu(\tilde {X_0}). $$
This formula is known under the name of L\^e-Greuel formula
\cite{Le}, \cite{Gr}.

Let us denote by  
$\phi_i(x) \in \Phi, 1 \leq i \leq \mu$
a basis of $\Phi.$ 

\begin{prop}We have the following free ${\mathcal O}_S$ module of rank $\mu,$
$$ \Phi(s)= 
\frac{ {\mathcal O}_{{\tilde X} \times S}}{{\mathcal O}_{{\tilde X} 
\times S}
\langle F_2(x,t),\cdots, F_k(x,t) \rangle + I_{C_0}(t)}.
$$
\label{prop50}
\end{prop}
{\bf Proof}
We reproduce the argument by \cite{bru}, Lemma A 1.
First of all we see that the module  $ \Phi(s)$ is a finitely generated
${\mathcal O}_S$ module. This can be shown by 
a combination of the Weierstra\ss-Malgrange preparation theorem and the
fact that for each fixed $s \in S$ the space 
$$ \frac{ {\mathcal O}_{\tilde X}}
{{\mathcal O}_{{\tilde X}}
\langle F_2(x,t),\cdots, F_k(x,t) \rangle + I_{C_0}(t)},
\leqno(5.3)$$
is a finite dimensional ($\leq \mu$) $\bC$ vector space (see \cite{Sz1}).

The above space $(5.3)$ is isomorphic to 
the direct sum of $\bC$ vector spaces,
$$ \bigoplus_{ \{x'; (x',s) \in Cr(\vec F)\} } 
\frac{{\mathcal O}_{\tilde X, x'}}{{\mathcal O}_{\tilde X, x'}
\langle F_2(x,t),\cdots, F_k(x,t) \rangle_{x'} 
+ I_{C_0}(t)}_{x'}.$$
Since this direct sum has dimension $\mu=$ 
the multiplicity of the critical point
$(x,u)=0$ of the height function on  $X_0,$ as mentioned at the very beginning of the paper,
it follows that 
$\{\phi_i(x)\}_{0 \leq i \leq \tau}$ form in fact a  
$\bC$ basis of $(5.3).$ Now we see that they form in fact
$\Phi(s)$ freely. If not, there exist holomorphic 
functions $\{a_i(s)\}_{0 \leq i \leq \tau}$ such that
$\sum^\tau_{i=0} a_i(s) \phi_i(x) =0$ in $\Phi(s).$ 
It would contradict the fact that for each fixed $s$,  
$\{\phi_i(x)\}_{0 \leq i \leq \tau}$ are linearly independent in
$(5.3).$ {\bf Q.E.D.}

\par
Let us consider the multiplication table
$$(F_1(x,t)-u) \phi_i(x\equiv \sum_{\ell=0}^\tau \rho_i^\ell(s) \phi_\ell(x)
\; mod ({\mathcal O}_{{\tilde X} \times S}
\langle F_2(x,t),\cdots, F_k(x,t) \rangle + I_{C_0}(t)).
\leqno(5.4)$$
Thus the matrix
$$ P(s) := (\rho_i^\ell(s))_{0 \leq i,\ell \leq \tau} =
(\tilde \rho_i^\ell(t) - u \cdot \delta_{i,\ell})_{0 \leq i,\ell \leq \tau},
\leqno(5.5)$$
is defined.
In analogy with $(3.3)$, we define another multiplication table
$$ \phi_i(x) \phi_j(x) \equiv \sum_{\ell=0}^\tau w^\ell_{i,j}(t)  
\phi_\ell(x)
\;\;mod ({\mathcal O}_{{\tilde X} \times S}
\langle F_2(x,t),\cdots, F_k(x,t) \rangle + I_{C_0}(t)). \leqno(5.6)$$
We will denote by $W^c(t)$ the matrix  
$(w^c_{\ell,b}(t))_{0 \leq \ell,b \leq \tau}.$
Hence,
$$(F_1(x,t)-u) \phi_a(x) \phi_b(x)
\equiv \sum_{\ell=0}^\tau \rho_a^\ell(s) \phi_\ell(x) \phi_b(x)  \leqno(5.7)$$
$$ \equiv \sum_{\ell=0}^\tau \rho_a^\ell(s) \sum_{c=0}^\tau 
w^c_{\ell,b}(t) \phi_c(x) \;\;mod ({\mathcal O}_{{\tilde X} \times S}
\langle F_2(x,t),\cdots, F_k(x,t) \rangle + I_{C_0}(t)). \leqno(5.7)$$
$$ \zeta_c(t) := tr(\phi_c(x) \cdot) = \sum_{\ell=0}^\tau w^\ell_{c,\ell}(t). 
\leqno(5.8)$$
Thus
$$tr ((F_1(x,t)-u) \phi_a(x) \phi_b(x) \cdot)
=  \sum_{\ell=0}^\tau \rho_a^\ell(s) \sum_{c=0}^\tau 
w^c_{\ell,b}(t) \zeta_c(t). \leqno(5.9)$$
We introduce the notation,
$$T(t)=\sum^\tau_{c=0} \zeta_c(t) W^c(t). \leqno(5.10)$$
From here on we assume that
$|s|$ is small enough and  denote by $\tilde X
=\{x \in \bC^n; |x| \leq  \delta\}$ a closed ball 
such that all critical points of 
$F_1(x,t)-u$ on $F_2(x,t)=\cdots= F_k(x,t)=0$ 
are located inside $\tilde X.$

In combining the results of \cite{Sz1}, Thoerem 2.1, Theorem 3.1,
with our above arguments we get the following.
\begin{thm}

1. The discriminant set of the deformation of projection
$X_t$ is given by the matrix $(5.5),$
$$ D=\{s \in S; det(P(s))=0\}.  \leqno(5.11)$$

2. $\{$ number of positive critical points of
$F_1(x,t)-u$ on $F_2(x,t)=\cdots= F_k(x,t)=0$, 
$x \in \tilde X \cap \bR^n$$\}$
- $\{$ number of negative critical points of
$F_1(x,t)-u$ on $F_2(x,t)=\cdots= F_k(x,t)=0$, 
$x \in \tilde X \cap \bR^n$ $\}$
$$= sign(P(s) \cdot T(t)).$$
\label{thm51}
\end{thm}

In opposition to the case $k=1$, we cannot write down a simple
formula for Euler characterisic of closures of semi-algebraic sets,
$$W_\ast =\{x \in \tilde X \cap \bR^n;
F_1(x,t)-u \ast 0 , F_2(x,t)=\cdots= F_k(x,t)=0\},$$
with $\ast = \geq, \leq , =$. As a matter of fact, it is quite easy
to establish an analogous theorem to \cite{Sz1} Theorem 3.3
on $\chi(W_{\geq 0}) \pm \chi(W_{\leq 0})$
by the  aid of matrices introduced above. We leave this task as an exercise
in view of complicated form of the analogy to the Hessian.

The bifurcation set $B_{F_1}$
is defined as
$B_{F_1} := \{ t \in T$; number of critical points of
$F_1(x,t)-u$ on $F_2(x,t)=\cdots= F_k(x,t)=0$ is strictly less than $\mu \} 
\setminus B_M.$
Here $B_M$ denotes the Maxwell set of $F_1(x,t)-u$, namely
$B_{M} := \{ t \in T$; two critical values of
$F_1(x,t)-u$ on $F_2(x,t)=\cdots= F_k(x,t)=0$ coincides $\}.$
\begin{prop}
The bifurcation set has the following expression
$$B_{F_1} = \{t \in T; det\; T(t)=0 \}. \leqno(5.12)$$
\label{prop52}
\end{prop}
{\bf Proof}

We consider the critical set
$$ C_0(t):= \{x \in \tilde X  ; dF_1(x,t) \wedge dF_2(x,t) \wedge \cdots \wedge dF_k(x,t)=0,
F_2(x,t)= \cdots =F_k(x,t)=0\}.$$ Here we remark that the critical set $C_0(t)$
has codimension $n$ in $\tilde X$ for a fixed generic value $t$ and it is a
set of points.
After \cite{Sz1} Corollary 2.5,
the rank of $T(t)$ is equal to the number of points $ \{p \in  C_0(t)\}.$
Therefore $T(t)$ degenerates if and only if $| C_0(t)| <\mu$
which means our statement.
{\bf Q.E.D.}

Regretfully, to the moment we cannot state how to deduce
the basis of $Der_S(log\;D)$ from the matrix $P(s).$
Consequently we cannot establish the relationship between 
the Gauss-Manin system and the topology of the real
algebraic sets.  This fact is due to the situation mentioned in the
Remark ~\ref{rmk54} below.

To remedy the situation, we state a
proposition on the multiplication table and the coefficients
to the Gauss-Manin system.

Let us consider the multiplication between $\phi_i$ and $\vec v_j$ by
the following way,
$$ \frac {\partial (\phi_i(x) h_{j,p}(x,s))}{\partial x_p}\equiv 
\sum^\tau_{r=0} R^r_{i,j}(s) \phi_r(x) mod ({\mathcal O}_{{\tilde X} \times S}
\langle F_2(x,t),\cdots, F_k(x,t) \rangle + I_{C_0}(t)).\leqno(5.13)$$
Here $ \check{\vec v_j}= \sum_{p=1}^n
h_{j,p}(x,s)\frac {\partial }{\partial x_p}$ denotes the vector field that
has been defined  in Lemma~\ref{lem32}.

\begin{prop}
The Gauss-Manin system for the period integrals 
$I_{\phi_i}^{(-\b1)}(s)$ introduced in the Theorem ~\ref{thm3}
is expressed by means of multiplication tables $(5.6)$
and $(5.13)$ as follows,
$$\vec v_j(I_{\phi_i}^{(-\b1)}(s)) = \sum^\mu_{\ell =1} 
\left((tr\;M^0)\cdot w^\ell_{i,j}(s)+ 
R^\ell_{i,j}(s) \right) I_{\phi_\ell}^{(-\b1)}(s)
\;\; 1 \leq j,q \leq \mu. $$
Here  $tr\;M^0$ stands for the trace of the non-degenerate matrix $M^0$
defined in Lemma ~\ref{lem31}.
\label{prop53}
\end{prop}

{\bf Proof}

First of all we remark the following chain of equalities,
$$\vec v_j \left( \int_{t(\gamma)} \phi_i(x)\vec F^{-\b1} dx \right)
= \int_{t(\gamma)} \phi_i(x) \left(\sum_{\ell=1}^k\sigma_j^\ell(s)
\frac{\partial }{\partial s_\ell} \vec F^{-\b1}\right) dx$$
$$ =\sum_{q=1}^k\ \int_{t(\gamma)} \phi_i(x)  F_q^{-1}\vec F^{-\b1}(\sum_{\ell=1}^k \sigma_j^\ell(s)
 \frac{\partial F_q(x,s)}{\partial s_\ell}) 
dx.$$
Here we remember Lemmata ~\ref{lem31}, ~\ref{lem32} and see that the above expression equals to
$$ \int_{t(\gamma)} \phi_i(x)  \phi_j(x) \left(
\sum_{q=1}^k F_q^{-1} (\sum_{\ell=1}^k C_q^\ell F_\ell(x,s) 
+ m_q(x, \vec F(x,s)) )\right) \vec F^{-\b1} dx$$
$$- \int_{t(\gamma)} \phi_i(x) \sum_{q=1}^k (-1)^{q-1} dF_q \wedge i_{\check{\vec v_j}}(dx)
 F_q^{-1}\vec F^{-\b1}.$$
As the terms with $C_q^\ell$, $\ell \not = q$ $($resp. terms with
$ m_q(x, \vec F(x,s))$ $\in$ ${\mathcal O}_{\tilde X}\otimes m^2_{S}$ )
vanish because of the lack of residues
along $F_\ell(x,s) =0$ (resp. some other $F_r(x,s) =0$), the last expression in its turn equals to
$$\int_{t(\gamma)} \phi_i(x) \phi_j(x)(\sum_{q=1}^k C_q^q F_q(x,s)) 
F_q^{-1}\vec F^{-\b1}  dx  + \int_{t(\gamma)} d(\phi_i(x) i_{\check{\vec v_j}}(dx))F^{-\b1} $$
$$ = (\sum_{q=1}^k C_q^q)\int_{t(\gamma)} \phi_i(x) \phi_j(x)\vec F^{-\b1} dx
+ \sum_{\ell=1}^k \int_{t(\gamma)} R^\ell_{i,j}(s) \phi_\ell(x) \vec F^{-\b1} dx$$
$$ = \sum^\mu_{\ell =1} \left((tr\;M^0)\cdot w^\ell_{i,j}(t)
+ R^\ell_{i,j}(s) \right) I_{\phi_\ell}^{(-\b1)}(s).$$

{\bf Q.E.D.}

\begin{remark}
{\rm The rank of $\bC-$module  of Leray coboundaries
 $t(\gamma) \in$
$H_{n}(\tilde X \setminus
\cup_{q=1}^k \{x \in \tilde X; F_q(x,s)=0 \})$ is equal to 
$\mu(\tilde X_0)$: the Milnor number of the singularity $\tilde X_0$
due to the tube operation isomorphism $t:$ defined in Lemma ~\ref{lem26}.
In view of the L\^e-Greuel formula mentioned in connection with $(5.2)$,  
the dimension $\mu$
of the space $\Phi$ is bigger than $\mu(\tilde X_0)$
as it represents the sum of the ranks of 
$(n-k)-$ dimensional cycles and $(n-k+1)-$dimensional cycles. Thus
we have no exact duality between the integrands and the integration cycles.
This means that the Gauss-Manin system of the above Proposition
~\ref{prop53} is defined only for the Riemann period matrix
of size $\mu \times \mu(\tilde X_0)$. 

To get the the Gauss-Manin system defined for the Riemann period matrix
of size $\mu(\tilde X_0) \times \mu(\tilde X_0)$, one need to consider
the multiplication table on the Brieskorn-Greuel lattice
$$ {\mathcal H}'':= \frac{\Omega_{\tilde X}^n}{dF_1(x,s)
\wedge \cdots \wedge dF_k(x,s)\wedge d\Omega_{\tilde X}^{n-k-1} + 
 \langle F_1(x,s),\cdots, F_k(x,s) \rangle \Omega_{\tilde X}^n},$$ 
that is known to be a ${\mathcal O}_S$ free module of
rank $\mu(\tilde X_0).$ This procedure can be done
in an analogous way to that in Proposition ~\ref{prop53}.
For the case of quasihomogeneous i.c.i.s., the concrete calculus
of the the Gauss-Manin system is done 
by means of Brieskorn-Greuel lattice in \cite{Tan99}.}
\label{rmk54}
\end{remark}

{ \center{\section{Examples}} }
1. Let us consider the simplest example
of the Pham-Brieskorn singularity,
$$F(x_1,x_2) = x_1^3 +x_2^3 +u + bx_1x_2+ c x_1 +d x_2, $$
with deformation parameters $s=(u,t) = (u,b,c,d).$
We calculate the data $(4.4)$, $(2.4)$, $(4.10)$,
$(4.11)$ as follows.
$$\tau^1=\left[ \begin {array}{cccc} 
1&0&0&0\\\noalign{\medskip}0&-1/3\,d&0&1/9\,bc
\\\noalign{\medskip}0&0&-1/3\,c&1/9\,bd\\\noalign{\medskip}0&1/9\,bc&1/9\,bd&1/9
\,dc\end {array} \right]$$

$$
\tau^2=
 \left[ \begin {array}{cccc} 0&1&0&0\\\noalign{\medskip}1&0&0&1/9\,{b}^{2}
\\\noalign{\medskip}0&0&-1/3\,b
&-1/3\,c\\\noalign{\medskip}0&1/9\,{b}^{2}&-1/3\,
c&1/9\,bd\end {array} \right]
$$

$$
\tau^3=\left[ \begin {array}{cccc} 0&0&1&0\\\noalign{\medskip}0&-1/3\,b&0&-1/3\,d
\\\noalign{\medskip}1&0&0&1/9\,{b}^{2}\\\noalign{\medskip}0&-1/3\,d&1/9\,{b}^{2}
&1/9\,bc\end {array} \right]
$$

$$
\tau^4=
 \left[ \begin {array}{cccc} 0&0&0&1\\\noalign{\medskip}0&0&1&0
\\\noalign{\medskip}0&1&0&0\\\noalign{\medskip}1&0&0&1/9\,{b}^{2}\end {array}
 \right]$$

$$
 \Sigma(s)=\left[ \begin {array}{cccc} 3\,u&2\,d&2\,c&b\\\noalign{\medskip}-2/3\,{d}^{2}+1
/9\,{b}^{2}c&3\,u+1/9\,{b}^{3}&-bd&2\,c\\\noalign{\medskip}-2/3\,{c}^{2}+1/9\,{b
}^{2}d&-bc&3\,u+1/9\,{b}^{3}&2\,d\\\noalign{\medskip}5/9\,bcd&-2/3\,{c}^{2}+1/3
\,{b}^{2}d&-2/3\,{d}^{2}+1/3\,{b}^{2}c&3\,u+1/9\,{b}^{3}\end {array} \right]
$$
$$
B^H=
 \left[ \begin {array}{cccc} 8\,{b}^{2}&16\,bc&16\,bd&{b}^{4}+16\,dc
\\\noalign{\medskip}16\,bc&-8\,{b}^{2}d&{b}^{4}+16\,dc&8/3\,{b}^{3}c-
16/3\,b{d}^{2}\\\noalign{\medskip}16\,bd&{b}^{4}+16\,dc&-8\,{b}^{2}c&8
/3\,{b}^{3}d-16/3\,b{c}^{2}\\\noalign{\medskip}{b}^{4}+16\,dc&8/3\,{b}
^{3}c-16/3\,b{d}^{2}&8/3\,{b}^{3}d-16/3\,b{c}^{2}&{\frac {56}{9}}\,{b}
^{2}dc+1/9\,{b}^{6})\end {array}\right]$$

$$ B^{HF}(s) =\Sigma(s) \cdot B^H=$$
$${\tiny \left[ \begin {array}{cc} (1.1):
24\,u{b}^{2}+80\,bcd+{b}^{5}&(1.2):-{\frac {64}
{3}}\,{d}^{2}{b}^{2}+14/3\,{b}^{4}c+48\,ubc+32\,d{c}^{2}\\
(1.3):-{\frac {64}{
3}}\,{c}^{2}{b}^{2}+14/3\,{b}^{4}d+48\,ubd+32\,{d}^{2}c&
(1.4):{\frac {152}{9
}}\,{b}^{3}dc-{\frac {32}{3}}\,b{c}^{3}-{\frac {32}{3}}\,b{d}^{3}+3\,u
{b}^{4}+48\,udc+1/9\,{b}^{7}\\
(2.1):{\medskip}-{\frac {64}{3}}\,{d}^
{2}{b}^{2}+14/3\,{b}^{4}c+48\,ubc+32\,d{c}^{2}&
(2.2):-{\frac {112}{3}}\,bc{d
}^{2}+{\frac {64}{9}}\,{b}^{3}{c}^{2}-24\,{b}^{2}du-{\frac {17}{9}}\,{
b}^{5}d\\
(2.3):{\frac {152}{9}}\,{b}^{3}dc-{\frac {32}{3}}\,b{c}^{3}-{\frac {
32}{3}}\,b{d}^{3}+3\,u{b}^{4}+48\,udc+1/9\,{b}^{7}&
(2.4)
\\
(3.1):{\medskip}
-{\frac {64}{3}}\,{c}^{2}{b}^{2}+14/3\,{b}^{4}d+48
\,ubd+32\,{d}^{2}c&
(3.2):{\frac {152}{9}}\,{b}^{3}dc-{\frac {32}{3}}\,b{c}^{
3}-{\frac {32}{3}}\,b{d}^{3}+3\,u{b}^{4}+48\,udc+1/9\,{b}^{7}\\
(3.3):
-{\frac{112}{3}}\,bd{c}^{2}+{\frac {64}{9}}\,{b}^{3}{d}^{2}-{\frac {17}{9}}\,
{b}^{5}c-24\,{b}^{2}cu&
(3.4)
\\
(4.1):{\medskip}{\frac {152}{9}}\,{
b}^{3}dc-{\frac {32}{3}}\,b{c}^{3}-{\frac {32}{3}}\,b{d}^{3}+3\,u{b}^{
4}+48\,udc+1/9\,{b}^{7}&
(4.2)\\
(4.3)&
(4.4)\end {array} \right]},$$
where
$$(2.4) = (4.2)= (3.4) =(4.3)$$
$$=-{\frac {106}{27}}\,{c}^{2}{b}^{4}-{
\frac {32}{3}}\,{c}^{3}d+{\frac {17}{27}}\,{b}^{6}d+{\frac {176}{9}}\,
{b}^{2}c{d}^{2}+8\,u{b}^{3}d-16\,ub{c}^{2},$$
$$ (4.4)={\frac {245}{81}}\,{b}^{5}c
d+16\,b{c}^{2}{d}^{2}-{\frac {32}{9}}\,{c}^{3}{b}^{3}-{\frac {32}{9}}
{b}^{3}{d}^{3}+{\frac {56}{3}}\,u{b}^{2}dc+1/3\,u{b}^{6}+{\frac {1}{
81}}\,{b}^{9}.$$
After Theorem ~\ref{thm42} the signature of this matrix gives us
the Euler characteristic of real algebraic sets
defined by $F(x,s) \geq , \leq , = 0.$

We calculate the determinants of these matrices.
$$
det(B^H)=1/9\, 
\left(256\,{b}^{2}{d}^{3}+768\,{d}^{2}{c}^{2}+96\,{b}^{4}dc-{b}^{8}+256\,
{c}^{3}{b}^{2} \right) ^{2},
$$
$$
det(\Sigma (s))=
8/3\,{b}^{2}{c}^{4}d-{\frac {1}{243}}\,{b}^{8}cd+8/3\,{d}^{4}c{b}^{2}+{\frac {23}{27}}\,{
b}^{4}{d}^{2}{c}^{2}+32\,ub{c}^{2}{d}^{2}-{\frac {11}{9}}\,u{b}^{5}cd-30\,{u}^{2}{b}^{2}d
c-{\frac {1}{243}}\,{b}^{6}{d}^{3}$$
$$-{\frac {1}{243}}\,{b}^{6}{c}^{3}-{\frac {32}{9}}\,{d}^
{3}{c}^{3}+24\,{u}^{2}{d}^{3}+1/3\,{u}^{2}{b}^{6}+9\,{u}^{3}{b}^{3}+{\frac {1}{243}}\,u{b
}^{9}-{\frac {20}{9}}\,u{c}^{3}{b}^{3}$$
$$-{\frac {20}{9}}\,u{b}^{3}{d}^{3}+24\,{c}^{3}{u}^{2
}+81\,{u}^{4}+{\frac {16}{9}}\,{d}^{6}+{\frac {16}{9}}\,{c}^{6}
$$

The discriminant of the polynomial $det(\Sigma)(s)$
with respect to the variable $u$ is calculated as follows,
$$Dscrim(det(\Sigma),u)=
27\, ( d-c )^{2}( 
{d}^{2}+dc+{c}^{2}) ^{2} ( 256\,{b}^{2}
{d}^{3}+768\,{d}^{2}{c}^{2}+96\,{b}^{4}dc-{b}^{8}+
256\,{c}^{3}{b}^{2})^{3}.
$$

These results combined with the Proposition~\ref{prop52}
calculate the Maxwell set,
$$  M= \{s \in \bC^3;   
\left( d-c \right) ^{2} \left( {d}^{2}+dc+{c}^{2} \right) ^{2}
=0 \}.$$

{\bf Example 2. The versal deformation of the singularity $E_6$.}

We consider the following deformation,
$$ F(x,y, t)+u= x^3+y^4 + gxy^2 + dy^2+ cxy + by + ax+u.$$
with $t=(a,b,c,d,g).$ As $F(x,y,0)$ is a quasihomogeneous polynomial in 
$(x,y)$,  we attribute to the deofromation parameters $(u,t) \in S$
corresponding quasihomogeneous weights. This means that there is a $\bC^\ast$
action on the space of deformation parameters $S.$
This allows us to consider $\tilde X= \bC^2$ , $S= \bC^6$
in the arguments of \S 4.

Thus we deal with the global parameter values $t \in \bC^5.$
Essentially all the informations on the multiplication table $(4.3)$
are contained in the following equivalence relations,

$x^2 \equiv
-1/3\,g{y}^{2}-1/3\,cy-1/3\,a 
\;\;mod (d_xF(x,y,t), d_y F(x,y,t))$

$y^3 \equiv( -1/2\,gy-1/4\,c ) x-1/2\,dy-1/4\,b$

$x^2y \equiv( 1/12\,gc+1/6\,{g}^{2}y ) x+1/12\,gb-1/3\,c{y}^{2}+ ( 1/6\,dg-1/3\,
a ) y$

$x^2y^2 \equiv
( 1/6\,{g}^{2}{y}^{2}+1/4\,gcy+1/12\,{c}^{2} ) x+ ( 1/6\,dg-1/3\,a
 ) {y}^{2}+ ( 1/12\,gb+1/6\,cd ) y+1/12\,cb$

$xy^3 \equiv
 (  ( -1/2\,d-1/12\,{g}^{3} ) y-1/24\,{g}^{2}c-1/4\,b ) x+1/4\,
gc{y}^{2}+ ( 1/6\,ga+1/12\,{c}^{2}-1/12\,{g}^{2}d ) y-1/24\,{g}^{2}b+1/12
\,ca$

$x^2y^3 \equiv
( 1/4\,gc{y}^{2}+ ( 1/3\,ga-1/3\,{g}^{2}d+1/12\,{c}^{2}-1/36\,{g}^{5}
 ) y-{\frac {1}{72}}\,{g}^{4}c+1/6\,ca-1/12\,{g}^{2}b-1/12\,dgc ) x+
 ( 1/12\,gb+1/6\,cd+1/12\,{g}^{3}c ) {y}^{2}+ ( -1/36\,{g}^{4}d+1/18
\,{g}^{3}a+1/3\,ad+1/36\,{g}^{2}{c}^{2}-1/6\,{d}^{2}g+1/12\,cb ) y+1/36\,{g}^{2
}ca-1/12\,dgb-{\frac {1}{72}}\,{g}^{4}b+1/6\,ab$

$xy^4 \equiv (  ( -1/12\,{g}^{3}-1/2\,d ) {y}^{2}+ ( -1/4\,b-1/6\,{g}^{2}c
 ) y-1/16\,g{c}^{2} ) x+ ( 1/12\,{c}^{2}-1/12\,{g}^{2}d+1/6\,ga
 ) {y}^{2}+ ( -1/24\,{g}^{2}b-1/8\,dgc+1/12\,ca ) y-1/16\,gcb$

$$x^2y^4 \equiv
(( 1/12\,{c}^{2}+1/3\,ga-1/36\,{g}^{5}-1/3\,{g}^{2}d ) {y}^{2}+$$
$$(-{\frac {11}{144}}\,{g}^{4}c-1/8\,{g}^{2}b+1/6\,ca-{\frac {7}{24}}\,dgc)y$$
$$-1/12\,gcb-1/32\,{g}^{3}{c}^{2}-1/24\,{c}^{2}d ) x$$
$$+ ( 1/3\,ad-1/
6\,{d}^{2}g+1/18\,{g}^{3}a-1/36\,{g}^{4}d+1/12\,cb+{\frac {13}{144}}\,{g}^{2}{c}^{2}
) {y}^{2}$$
$$+ ( -1/12\,c{d}^{2}-1/16\,{g}^{3}cd+1/6\,ab-1/8\,dgb+1/48\,g{c}^
{3}+{\frac {5}{72}}\,{g}^{2}ca-{\frac {1}{72}}\,{g}^{4}b) y$$
$$-1/48\,g{b}^{2}-1/
24\,cdb+1/48\,g{c}^{2}a-1/32\,{g}^{3}cb.$$
We can write down these results in the form of matrices $(4.4)$
and the polynomials $\zeta_k(t)$,
$1 \leq k \leq 6$, $(4.6)$,
$$
\tau^1=\left[ \begin {array}{cccccc}
  1&0&0&0& 0 & 0 \\
  0 & -\frac{a}{3} & 0 &    \frac{b\,g}{12} & 0 & \frac{b\,c}{12} \\
  0 & 0 & 0 & 0 & -\frac{b}{4} &
   \frac{a\,c}{12} - \frac{b\,g^2}{24} \\
 0 & \frac{b\,g}{12} & 0 & \frac{b\,c}{12} &
   \frac{a\,c}{12} - \frac{b\,g^2}{24} &p_1\\
  0 & 0 & -\frac{b}{4} & \frac{a\,c}{12} -
    \frac{b\,g^2}{24} & 0 &
   -\frac{  b\,c\,g  }{16} \\
  0 & \frac{b\,c}{12} &
   \frac{a\,c}{12} - \frac{b\,g^2}{24} &
   p_1 &
   -\frac{ b\,c\,g  }{16} &
   -\frac{b\,c\,d}{24} -
    \frac{b^2\,g}{48} + \frac{a\,c^2\,g}{48} -
    \frac{b\,c\,g^3}{48}  \end{array} \right]$$
where $ p_1=\frac{a\,b}{12} - \frac{b\,d\,g}{24} +
    \frac{a\,c\,g^2}{72} - \frac{b\,g^4}{144}$,
$$ \zeta_1(t) =6 .$$

$$
\tau^2=\left[ \begin {array}{cccccc}
  0 & 1 & 0 & 0 & 0 & 0 \\
 1 & 0 & 0 & \frac{c\,g}{12} & 0 & \frac{c^2}{12} \\
 0 & 0 & 0 & 0 & -\frac{c}{4} & -\frac{b}{4} -
\frac{c\,g^2}{24} \\ 0  &
  \frac{c\,g}{12} & 0 & \frac{c^2}{12} & -\frac{b}{4} - \frac{c\,g^2}{24} &
   p_2 \\
  0 & 0 & -\frac{c}{4} & -\frac{b}{4} - \frac{c\,g^2}{24} & 0 & -\frac{
c^2\,g  }{16} \\
  0 & \frac{c^2}{12} & -\frac{b}{4} - \frac{c\,g^2}{24} &
  p_2 &
   -\frac{ c^2\,g  }{16} &
   -\frac{c^2\,d }{24} - \frac{b\,c\,g}{12} - \frac{c^2\,g^3}{48}
\end {array} \right], $$
where $p_2=\frac{a\,c}{12} - \frac{c\,d\,g}{24} -
\frac{b\,g^2}{24} - \frac{c\,g^4}{144}$,
$$ \zeta_2=\frac{g^2}{3}.$$

$$\tau^3=
\left[ \begin {array}{cccccc} 
 0 & 0 & 1 & 0 & 0 & 0\\
 0 & -\frac{c}{3} & 0 & -\frac{a}{3} +
\frac{d\,g}{6} & 0 & \frac{c\,d}{6} + \frac{b\,g}{12}\\
 1 & 0 & 0 & 0 & -\frac{d}{2} & \frac{c^2}{12} + \frac{a\,g}{6} -
\frac{d\,g^2}{12} \\
  0 & -\frac{a}{3} + \frac{d\,g}{6} & 0 & \frac{c\,d}{6} + \frac{b\,g}{12} &
\frac{c^2}{12} + \frac{a\,g}{6}
- \frac{d\,g^2}{12} & p_3 \\
   0 & 0 & -\frac{d}{2} & \frac{c^2}{12} + \frac{a\,g}{6} -
\frac{d\,g^2}{12} & -\frac{b}{4} &
   \frac{a\,c}{12} - \frac{c\,d\,g}{8} - \frac{b\,g^2}{24}\\
  0 & \frac{c\,d}{6} + \frac{b\,g}{12} & \frac{c^2}{12} + \frac{a\,g}{6} -
\frac{d\,g^2}{12} & p_3 &
   \frac{a\,c}{12} - \frac{c\,d\,g}{8} - \frac{b\,g^2}{24} & q_3\end {array} \right],$$
where $$p_3=\frac{b\,c}{12} + \frac{a\,d}{6} - \frac{d^2\,g}{12} +
\frac{c^2\,g^2}{72} + \frac{a\,g^3}{36} - \frac{d\,g^4}{72},$$ $$
q_3=\frac{a\,b}{12} - \frac{c\,d^2}{12} + \frac{c^3\,g}{48} -
\frac{b\,d\,g}{12} + \frac{a\,c\,g^2}{18} - \frac{c\,d\,g^3}{24} -
  \frac{b\,g^4}{144}.$$
$$ \zeta_3(t) =0 .$$

$$
\tau^4=\left[ \begin {array}{cccccc}
 0 & 0 & 0 & 1 & 0 & 0 \\
  0 & 0 & 1 & \frac{g^2}{6} & 0 & \frac{c\,g}{4} \\
   0 & 1 & 0 & 0 & -\frac{g}{2} & -\frac{d}{2} -\frac{g^3}{12} \\
  1 & \frac{g^2}{6} & 0 & \frac{c\,g}{4} & -\frac{d}{2} - \frac{g^3}{12} &
   \frac{c^2}{12} + \frac{a\,g}{6} - \frac{d\,g^2}{6} - \frac{g^5}{72} \\
  0 & 0 & -\frac{g}{2} & -\frac{d}{2} - \frac{g^3}{12} & -\frac{c}{4} &
-\frac{b}{4} - \frac{c\,g^2}{6}\\
  0 & \frac{c\,g}{4} & -\frac{d}{2} - \frac{g^3}{12} & \frac{c^2}{12} +
\frac{a\,g}{6} - \frac{d\,g^2}{6} - \frac{g^5}{72} &
   -\frac{b}{4} - \frac{c\,g^2}{6} & \frac{a\,c}{12} - \frac{c\,d\,g}{4} -
\frac{b\,g^2}{12} - \frac{7\,c\,g^4}{144}\end {array} \right].$$
$$\zeta_4(t)=\frac{5\,c\,g}{6}.$$

$$\tau^5=\left[ \begin {array}{cccccc}
  0 & 0 & 0 & 0 & 1 & 0\\
 0 & -\frac{g}{3} & 0 & -\frac{c}{3} & 0 & -\frac{a}{3} +
\frac{d\,g}{6} \\
 0 & 0 & 1 & 0 & 0 & \frac{c\,g}{4}\\
  0 & -\frac{c}{3} & 0 & -\frac{a}{3} + \frac{d\,g}{6} & \frac{c\,g}{4} &
  \frac{c\,d}{6} + \frac{b\,g}{12} + \frac{c\,g^3}{24}\\
  1 & 0 & 0 & \frac{c\,g}{4} & -\frac{d}{2} & \frac{c^2}{12} +
\frac{a\,g}{6} - \frac{d\,g^2}{12}\\
   0 & -\frac{a}{3} + \frac{d\,g}{6} & \frac{c\,g}{4} & \frac{c\,d}{6} +
\frac{b\,g}{12} + \frac{c\,g^3}{24} &
   \frac{c^2}{12} + \frac{a\,g}{6} - \frac{d\,g^2}{12} & q_5 \end{array} \right], $$
where $$q_5=\frac{b\,c}{12} + \frac{a\,d}{6} - \frac{d^2\,g}{12} +
\frac{11\,c^2\,g^2}{144} + \frac{a\,g^3}{36} -
\frac{d\,g^4}{72}.$$
$$\zeta_5(t)=-2\,d - \frac{g^3}{6}.$$

$$\tau^6=\left[ \begin {array}{cccccc}
  0 & 0 & 0 & 0 & 0 & 1 \\
  0 & 0 & 0 & 0 & 1 & \frac{g^2}{6} \\
  0 & 0 & 0 & 1 & 0 & 0 \\
  0 & 0 & 1 & \frac{g^2}{6} & 0 & \frac{c\,g}{4} \\
 0 & 1 & 0 & 0 & -\frac{g}{2} & -\frac{d}{2} - \frac{g^3}{12}\\
   1 & \frac{g^2}{6} & 0 & \frac{c\,g}{4} & -\frac{d}{2} - \frac{g^3}{12} &
   \frac{c^2}{12} + \frac{a\,g}{6} - \frac{d\,g^2}{6} - \frac{g^5}{72}\end {array} \right]$$
$$\zeta_6(t)=\frac{5\,c^2}{12} + \frac{2\,a\,g}{3} - \frac{d\,g^2}{2}
- \frac{g^5}{36}.$$

Finally we get the matrix $(4.5)$ as follows.
$$ T(t) = 6 \tau^1(t) + \zeta_3(t) \tau^3(t)+ \zeta_4(t) \tau^4(t) + \zeta_5(t) \tau^5(t) +
\zeta_6(t) \tau^6(t)=$$  
{\tiny
$$=\left[ \begin {array}{cccccc}
 6 &
  \frac{g^2}{3} &
   0 &
    \frac{5\,c\,g}{6} &
     -2\,d - \frac{g^3}{6} &
     T_{1,6}
\\
  \frac{g^2}{3} &
  -2\,a + \frac{2\,d\,g}{3} + \frac{g^4}{18}&
   \frac{5\,c\,g}{6} &
    \frac{2\,c\,d}{3} + \frac{b\,g}{2} + \frac{2\,c\,g^3}{9} &
    T_{2,5}  
&
      T_{2,6}
\\
 0 &
  \frac{5\,c\,g}{6} &
   -2\,d - \frac{g^3}{6} &
T_{3,4}    
&
      -\frac{3\,b}{2} - \frac{c\,g^2}{2} &
      T_{3,6}
\\
 \frac{5\,c\,g}{6} &
  \frac{2\,c\,d}{3} + \frac{b\,g}{2} + \frac{2\,c\,g^3}{9} &
T_{4,3}
   
&
   T_{4,4} &
T_{4,5}  &
 T_{4,6}
\\
  -2\,d - \frac{g^3}{6} &
  T_{5,2} 
&
    -\frac{3\,b}{2} - \frac{c\,g^2}{2} &
    T_{5,4}
&
T_{5,5}
&
T_{5,6}
\\
T_{6,1}&
  T_{6,2} &
 T_{6,3} &
  T_{6,4} &
  T_{6,5}&
   T_{6,6}\end{array} \right],$$}
where 
{\tiny
$$ \begin {array}{c}
 T_{1,6}= T_{2,5}=T_{3,4}=T_{4,3}=
T_{5,2}=T_{6,1}= \frac{5\,c^2}{12} +
\frac{2\,a\,g}{3} - \frac{d\,g^2}{2} - \frac{g^5}{36} ,
\\
T_{2,6}=T_{6,2}=\frac{b\,c}{2} + \frac{2\,a\,d}{3} - \frac{d^2\,g}{3} +
\frac{11\,c^2\,g^2}{36} + \frac{a\,g^3}{6} - \frac{d\,g^4}{9} -
  \frac{g^7}{216},\\
T_{4,4}= \frac{b\,c}{2} + \frac{2\,a\,d}{3} - \frac{d^2\,g}{3} +
\frac{11\,c^2\,g^2}{36} +\frac{a\,g^3}{6} - \frac{d\,g^4}{9} -
  \frac{g^7}{216}.\\ 
T_{4,5}=T_{5,4}= \frac{a\,c}{2} - \frac{11\,c\,d\,g}{12} -  \frac{b\,g^2}{3} -
\frac{c\,g^4}{8}\\
T_{3,6}=
T_{6,3}= \frac{a\,c}{2} - \frac{11\,c\,d\,g}{12} - \frac{b\,g^2}{3} -
  \frac{c\,g^4}{8}\\
 T_{5,5}=
d^2 - \frac{5\,c^2\,g}{12} - \frac{a\,g^2}{3} + \frac{d\,g^3}{3} +
\frac{g^6}{72}\\
  T_{4,6}=T_{6,4}=\frac{a\,b}{2} - \frac{c\,d^2}{3} +
\frac{25\,c^3\,g}{144} - \frac{5\,b\,d\,g}{12} +
\frac{5\,a\,c\,g^2}{12} -
  \frac{7\,c\,d\,g^3}{18} - \frac{5\,b\,g^4}{72} -
  \frac{c\,g^6}{36},\\
 T_{5,6}= T_{6,5}= -\frac{3\,c^2\,d}{8} - \frac{7\,b\,c\,g}{12} -
  \frac{2\,a\,d\,g}{3} + \frac{5\,d^2\,g^2}{12} - \frac{5\,c^2\,g^3}{24} -
  \frac{a\,g^4}{12} + \frac{5\,d\,g^5}{72} + \frac{g^8}{432} ,\\
  T_{6,6}=\\
\frac{5\,c^4}{144} -
   \frac{5\,b\,c\,d}{12} - \frac{a\,d^2}{3} - \frac{b^2\,g}{8} +
\frac{23\,a\,c^2\,g}{72} + \frac{d^3\,g}{6} +
  \frac{a^2\,g^2}{9} - \frac{35\,c^2\,d\,g^2}{72} - \frac{17\,b\,c\,g^3}{72}
- \frac{5\,a\,d\,g^3}{18} +
  \frac{d^2\,g^4}{8} -
  \frac{59\,c^2\,g^5}{864} - \frac{a\,g^6}{54} + \frac{d\,g^7}{72} +
  \frac{g^{10}}{2592}.\\
\end{array}$$}
It is a conceptually easy exercise to calculate further
$B^H(s)$ and $B^{HF}(s)$ to establish correspondence
between parameter value $s=(a,b,c,d,g,u)$ and the Euler characteristic 
of a semi-algebraic set defined by $F(x,y,t)+u$.

For instance, for the values
$$ -0.6 \leq a \leq 1, (b,c,d,g,u)=(-0.4,0.1, 0.1, -0.1, -10),$$
we calculate with computer (Mathematica computation achieved by Galina Filipuk)
$\chi(W_{\geq 0})$ $=$ $\chi(W_{\leq 0})$ $=$ $0,$
while for the values
$$ -1 \leq a \leq -0.8, (b,c,d,g,u)=(-0.4,0.1, 0.1, -0.1, -10),$$
we have $\chi(W_{\geq 0})=1, \chi(W_{\leq 0})=-1.$

For the values
$$ -1 \leq a \leq -0.8, (b,c,d,g,u)=(-0.4,0.1, 0.1, -0.1, 8.5),$$
we have $\chi(W_{\geq 0})=-1, \chi(W_{\leq 0})=1,$
and 
$$ -0.6 \leq a \leq 1, (b,c,d,g,u)=(-0.4,0.1, 0.1, -0.1, 8.5),$$
we have $\chi(W_{\geq 0})=0, \chi(W_{\leq 0})=0.$

It is worthy noticing that the first two cases (resp. last two cases)
give us examples of topologically different  isotopy types of the real curve
for the same sign combination of coefficients
$(-,-,+,+,-,-)$ (resp. $(-,-,+,+,-,+)$). These examples show the cases 
that Viro's patchworking method could not distinguish.


\vspace{\fill}

\noindent

\begin{flushleft}
\begin{minipage}[t]{6.2cm}
  \begin{center}
{\footnotesize Indepent University of Moscow\\
Bol'shoj Vlasievskij pereulok 11,\\
 Moscow, 121002,\\
Russia\\
{\it E-mails}:  tanabe at mccme.ru, tanabe at mpim-bonn.mpg.de \\}
\end{center}
\end{minipage}
\end{flushleft}

\begin{thebibliography}{99}
\bibitem{Al}{\sc A.G.Aleksandrov},
{\it
Nonisolated Saito singularities,}.
Math.USSR Sbornik  {\bf 65} (1990), No.2, pp.561-574.

%
\bibitem{bru}{\sc J.W.Bruce},
{\it Functions on discriminants,}J.London Math.Soc.
{\bf 30} (1984), pp.551-567.


%
\bibitem{degre}{\sc I. de Gregorio},
{\it Deformations of functions and F-manifolds,}
to appear in Bull.London Math.Soc.


%
\bibitem{Gory1}{\sc V.V.Goryunov},
 {\it Projections and vector fields that are tangent to
the discriminant of a complete intersection,}
Funct. Anal. Appl.   {\bf 22} (1988), No.2, pp.104-113.
%
\bibitem{Gory2}{\sc V.V.Goryunov},
 {\it Vector fields and functions on the discriminants of complete intersections and bifurcation diagrams of projections,} J. Soviet Math. {\bf 52} (1990), no. 4, pp.3231--3245.
%
 \bibitem{Gr}{\sc G.-M.Greuel},
 {\it Der Gau\ss-Manin-Zusammenhang isolierter Singularit\"aten von
 vollst\"andigen Durchschnitten,} Math. Ann. {\bf 214}  no.3, (1975),
 pp. 235-266.
%
\bibitem{GSch}{\sc Ph.Griffiths} and {\sc W.Schmid},
{\it Recent develpments in Hodge theory}, in Discrete 
subgroups of Lie groups and Applications to Moduli,
Oxford University Press, 1973, pp.31-127.
%
\bibitem{HT}{\sc R.C.Hwa} and {\sc V.L.Teplitz},
{\it Homology and Feynman integrals},
W.A.Benjamin, Inc., 1966.


\bibitem{Le}{\sc L\^e D\~ung Tr\'ang},
{\it Calculation of Milnor number of isolated 
singularity of complete intersection,}.
Funkts. Anal. Appl.   {\bf 8} (1974), No.2, pp.45-52.


\bibitem{Looi}{\sc E.Looijenga},
 {\it Isolated singular points on
complete intersections,} London Math. Soc. Lect. Notes Ser.,
1984, No. 77, 200pp.


\bibitem{Mal}{\sc B.Malgrange},
 {\it
Ideals of differentiable functions,} Tata Institute of Fundamental Research Studies in Mathematics, No. 3 
Tata Institute of Fundamental Research, Bombay; Oxford University Press, London 1967, vii + 107ppthese.
%
\bibitem{Saito0}
{\sc K.Saito},
 {\it Quasihomogene isolierte Singularit\"aten 
von Hyperfl\"achen,} Invent. Math.,{\bf 14} (1971), pp.123-142.
\bibitem{Saito1}
{\sc K.Saito},
 {\it Theory of logarithmic differential 
forms and logarithmic vector fields ,}
J.Fac. Sci. Tokyo, Sec. I A Math.,{\bf 27} (1980), pp.265-291.
%
\bibitem{Saito2}
{\sc K.Saito},
 {\it Period mapping associated to a primitive form
,} Publ. RIMS Kyoto Univ.,{\bf 19} (1983), pp.1231-1264.
%
\bibitem{Saito3}
{\sc K.Saito},
 {\it On the periods of primitive integrals I,}
RIMS preprint 412, Kyoto Univ. 1982. 

%
\bibitem{Sz1}{\sc Z.Szafraniec},
 {\it On topological invariants of real analytic singularities,}.
Math. Proc. Camb. Phil.Soc., {\bf 130} (2001), pp.13-24.

%
 \bibitem{Tan99}{\sc   S.Tanab\'e},   {\it
Transform\'ee  de Mellin des int\'egrales- fibres
associ\'ees  aux   singularit\'es isol\'ees
 d'intersection compl\`ete quasihomog\`enes,} 
Compositio Math. {\bf 130}(2002), no.2,
pp. 119-160.
%
\bibitem{Terao}{\sc H.Terao},
 {\it The bifurcation set and logarithmic vector fields ,}
Math. Ann.,   {\bf 263} (1983), No.3, pp.313-321.


%
 \bibitem{Vas}{\sc V.A.Vassiliev}, {\it Ramified integrals, singularities and
 Lacunas,} Kluwer Academic Publishers, Dordrecht, 1995.
%
\bibitem{Viro}{\sc O.Ya.Viro}, {\it
Real plane algebraic curves: constructions with controlled topology}, 
Leningrad Math. J. {\bf 1} (1990), no. 5, 1059--1134

\end{thebibliography}
\end{document}